\documentclass[12pt,a4,epsf]{article}
\usepackage{amsfonts}
\usepackage{graphicx}
\usepackage{latexsym}
\usepackage{epsfig}
\usepackage{amsmath}
\usepackage{amssymb}
\usepackage{color}
\usepackage{amsthm}
\usepackage[dvipdfm,CJKbookmarks,bookmarksopen=true,colorlinks=true,
linkcolor=blue,citecolor=blue,pdfstartview=FitH,pdftitle=title,pdfauthor=lixx]{hyperref}

\newtheorem{theo}{\sc Theorem}[section]
\newtheorem{lemm}{\sc Lemma}[section]
\newtheorem{prop}{\sc Proposition}[section]

\newtheorem{rema}{\sc Remark}[section]
\renewcommand{\theequation}{{\arabic{section}.\arabic{equation}}}

\newcommand{\be}{\begin{equation}}
\newcommand{\ee}{\end{equation}}
\newcommand{\bea}{\begin{eqnarray}}
\newcommand{\eea}{\end{eqnarray}}
\newcommand{\beas}{\begin{eqnarray*}}
\newcommand{\eeas}{\end{eqnarray*}}

\makeatletter      
\@addtoreset{equation}{section}
\makeatother       

\long \def\@makecation#1#2{ \vskip 10 pt
\setbox\@tempboxa\hbox{#1:#2} \ifdim \wd\@tempboxa >\hsize
\unhbox\@tempboxa\par \else \hbox to\hsize{\hfil\box\@temboxa\hfil}
\fi}

\def \be{\mbox{\boldmath $\beta$}}

\begin{document}

\title{\bf  Bounds smaller than the Fisher information for generalized linear models}

\author{
\  Lixing Zhu\footnote{The corresponding author.
Email: lzhu@hkbu.edu.hk.    Lixing  ZHU was supported by a grant
from the Research Grants Council of Hong Kong, and  a FRG grant from
Hong
Kong Baptist University, Hong Kong. In this research,
the methodology was developed and material organization of the manuscript
 was done by the first author. This manuscript is an updated version of the first manuscript that was
 submitted to arXiv: submit/0041278 on May 15, 2010. A special thank goes to
 Dr. Bing Li who has had intensive and stimulative discussions
 with us so that we have a  deeper thinking on possible super-efficiency.}\, \, and Zhenghui Feng   \\ 
\begin{tabular}{l}
{\small\it  Department of Mathematics, Hong Kong Baptist
University, Hong Kong, China}
\end{tabular}}
\date{}
\maketitle

\begin{abstract}
In this paper,  we propose
a parameter space augmentation approach  that is based on ``intentionally" introducing a pseudo-nuisance
parameter into generalized linear models 
for the purpose of variance reduction.
We first consider the parameter whose norm is equal to one. By  introducing a pseudo-nuisance parameter into models to be estimated, an extra estimation is asymptotically normal and is, more importantly, non-positively correlated to the estimation that asymptotically achieves the Fisher/quasi Fisher information. As such, the resulting estimation is asymptotically with smaller variance-covariance matrices than the Fisher/quasi Fisher information. For general cases where the norm of the parameter is not necessarily equal to one, two-stage quasi-likelihood procedures separately estimating the scalar and direction of the parameter are proposed. The traces of the limiting
variance-covariance matrices  are in general smaller
than or equal to that of the Fisher/quasi-Fisher information. 
We also
discuss the pros and cons of the new methodology, and possible extensions. As this methodology of parameter space augmentation is general, and then may be readily extended to handle, say, cluster data and correlated data, and other models.
\end{abstract}



\section{Introduction}\label{sec1}

For parametric regression models,  efficiency of
estimation is used as a criterion for estimation
accuracy and is a well investigated issue,
and relevant theories have been well developed
in the literature. For instance, for linear models, when
 the design is fixed,
  Gauss-Markov theorem shows that  the ordinary least squares estimation is the best linear unbiased estimation (BLUE, Markoff 1912).
  Plackett (1949) had a nice description on the
  contributions by Gauss, Laplace, Markov and others.
  Here `Best' means minimum variance. When the design
  is random, its limiting variance achieves the Fisher
  information. However, the unbiasedness of estimation
  is usually not achievable, particularly for generalized
   linear models (GLMs), a notion about `best' that can
    be explained for asymptotically unbiased estimation
     is also the Fisher information.

The Fisher information is a way of measuring the amount of
information that an observable random variable  carries
 about an unknown parameter ¦È upon which the probability
 of the observable random variable depends. Fisher (1922, p.316) presented the ``Criterion
 of Efficiency" that refers to large-sample behavior: ``when the
 distribution of the statistic tends to normality,
 that statistic is to be chosen which has the least
 probable error." ``Efficiency" is efficiency in
 summarizing information  in the large-sample
  case.


The class of generalized linear models (GLMs) is one of the most
widely used classes of regression models for statistical analyses. The
relevant theoretical investigations are rather complete and the
results of parameter estimation have been the standard results in
the literature and textbooks. As is well known, the GLM is of a
structure as
\begin{eqnarray}\label{glm1}
h(\mu)=\beta_0^TX
 \end{eqnarray}
 where $\mu=E(Y|X)$, $Y$ is the response and $X$ is the predictor
 vector of $p$ dimension, $\beta$ is the parameter of interest, and
 $h$ is a given link function. It is also often written as
\begin{eqnarray}\label{glm2}
Y=g(\beta_0^TX)+\varepsilon,
 \end{eqnarray}
with $E(\varepsilon|X)=0$ and a given function $g$ although
(\ref{glm1}) and (\ref{glm2})  are not exactly equivalent unless $g$
is invertible. Without loss of generality, assume that $E(X)=0$
throughout the present paper.

Suppose that an independent identically distributed sample $\{(x_1,
y_1), \cdots, (x_n, y_n)\}$ is available. To estimate the parameter
$\beta$, there are several standard methods and among them, the weighted
least squares and  the quasi-likelihood proposed by Wedderburn
(1974) are the most popularly used methods when the distribution of error is not parametric. The estimation of
$\beta$ via the quasi-likelihood can achieve the asymptotic efficiency: when the underlying distribution of error is natural exponential, the limiting variance-covariance matrix can be, or be proportional to, the Fisher information that attains the
Cram\'er-Rao lower bound (Cram\'er, 1946; Rao 1945). In effect, it is asymptotically
efficient among  all consistent solutions to linear estimating
equations, see Wedderburn (1974), McCullagh (1983), Jarrett (1984),
McLeish (1984), Firth (1987) and Godambe and Heyde (1987). A relevant reference is
Godambe (1960) who proved the optimum property of the Fisher
information in the likelihood estimation. When we do not impose, other than some moments of the error, specific distributional conditions, the limiting variance-covariance matrix is called
the quasi-Fisher information,  see e.g. Wefelmeyer (1996). These results have  long been benchmarks
to be compared for estimation efficiency of any new method. Le Cam (1986) is a comprehensive reference book about rigorous  proofs of the asymptotic optimality in this area. In
recent years, some approaches were developed to improve generalized
estimating equations by using quadratic inference functions for
cluster data, see Li and Lindsay (1996) and Qu, Lindsay and Li
(2000), and Wang and Qu (2009).

Note that all of existing estimation methods are within parametric
framework.   Nevertheless, this is a natural and seemingly necessary
way for a parametric problem. As such, it seems
 that we may not need to bother ourself to estimate anything else other
than the parameter of interest if there are no other nuisance
parameters in the model. From a heuristic perspective in both
practice and theory, if we did so, some extra estimation errors
would be created and estimation variation would be even larger.

However, our study shows that this is not always true. In this
paper, we propose a two-stage estimation method to estimate the
scalar $\|\beta_0\|$ and the direction $\beta_0/\|\beta_0\|$
separately to form a final estimator of $\beta_0$ where $\|\cdot\|$
is the Euclidean norm. Although both the stages are based on the
quasi-likelihood, the estimation for the direction
$\beta_0/\|\beta_0\|$ is a semiparametric quasi-likelihood by introducing a pseudo-nuisance parameter into models to be
estimated and
regarding the link function $g$ as its true value. In  the next section, we will give a detailed motivation to describe why we do this. The results are interesting: the traces of the limiting
variance-covariance matrices are smaller than or
equal to  those of the Fisher/quasi-Fisher information, and in some cases, particularly for $\beta_0$ with norm one, the matrices
themselves are smaller than or equal to the corresponding Fisher/quasi-Fisher information.
Also the reduction of variance is highly related to the correlations  among the
predictors. We will see this in the theorems in Section~\ref{sec4}.  Two relevant references are  Wang,
Xue, Zhu and Chong (2010) and Chang, Xue and Zhu (2010)  about the least squares.
Its special case is the linear model. Our result shows that although our estimation is only asymptotically unbiased, trace of
the variance-covariance matrix is smaller than what Gauss-Markov theorem provided
in independent identically distributed case.
In the next section, we shall present the
details of the estimation procedures and the results. Of course, as a trade-off,
to facilitate the use of nonparametric estimation, technical conditions on the predictor vector
are much more stronger than those the classical methods need.

It is worth pointing out that this parameter-space augmentation approach has a root from the estimation problems with nuisance parameters. For example, Piece (1982) pointed out that when  nuisance parameters in models are estimated, the estimation of the parameter of interest may be more efficient than that when  nuisance parameters are  regarded as given. The difference from his result is that as a methodology, we create a nuisance parameter to be estimated and regard the given link as its true value. This methodology may be useful for other estimation problems. We will have discussions in Section~\ref{sec3}.

The paper is organized as follows. In Section~\ref{sec2}, we shall in detail
describe the motivation of our method. Section~\ref{two-pahse} presents the  estimation procedures.
Section~\ref{sec4} presents the results. As the results are related to given
variance function of error and more regularity conditions than the
classical quasi-likelihood needs for technical purpose, we shall in Section~\ref{sec3}
discuss pros and cons of the new methodology and some extensions.
The proofs of the main results are relegated to Section~\ref{sec5},
the Appendix.

\section{Motivation of methodology development }\label{sec2}
\subsection{A brief review of the quasi-likelihood}
For the model of (\ref{glm2}), suppose that a sample $\{(x_1, y_1), \cdots,
(x_n, y_n)\}$ is available. In a simple case with a given variance function $v(X)=var(Y|X)$, the quasi-likelihood estimator $\hat \beta$ is the solution to the following equation
\begin{eqnarray}\label{gee1}
G(\beta)=\sum_{i=1}^n(y_i-g(\beta^Tx_i))g'(\beta^Tx_i)x_i/v(x_i),
 \end{eqnarray}
over all $\beta$. It is worth pointing out  that for this $v(\cdot)$, the above is equivalent to weighted least squares (WLS). When the variance can be written as $\sigma^2 v(X)$ where $\sigma^2$ is an unknown dispersion parameter, we can still use $v(X)$ in the estimating equation. However, when $v(X)$
is of the structure $v(\beta_0^TX)$ or $v(g(\beta_0^TX))$, it is different from the WLS, and we will use it in $G(\beta)$ in
the lieu of $v(X)$. As the proof is almost identical to the case
with $v(X)$, and the results are also almost unchanged, even simpler in form, for the
limiting variance of the estimator $\hat \beta$. Further, when function $v$ of $v(g(\beta_0^TX))$ is also unknown, we need to use a nonparametric estimation to replace it. The result is also similar.   The extensions will be discussed in  Section~\ref{sec3}. Thus, we for simplicity of presentation
use $v(X)$ throughout this paper.

From the projection theorem (Small
and Mcleish, 1994, p. 79), the quasi-likelihood equation $G(\beta)$
of (\ref{gee1}) is the optimal linear combination of
$(y_i-g(\beta^Tx_i))$'s. The asymptotic normality of $\hat \beta$
can be derived by using Taylor expansion and Slutsky theorem to
obtain an asymptotic presentation as
\begin{eqnarray}\label{gee2}
&&\frac1{\sqrt
n}\sum_{i=1}^n\varepsilon_ig'(\beta_0^Tx_i)x_i/v(x_i)\nonumber\\&=&\Big
(\frac 1n\sum_{i=1}^n(g'(\beta_0^Tx_i))^2/v(x_i)x_ix_i^T\Big )\sqrt
n(\hat
\beta-\beta_0)+o_p(1)\nonumber\\
&=:&V\sqrt n(\hat \beta-\beta_0)+o_p(1).
 \end{eqnarray}
with $V=E\Big(( g'(\beta_0^TX))^2/v(X)XX^T)$ and then
\begin{eqnarray}\label{gee3}
\sqrt n (\hat \beta-\beta_0)\Longrightarrow N(0, V^{-1})
 \end{eqnarray}
where ``$\Longrightarrow$" stands for convergence in distribution, and
 $V=E\Big( (g'(\beta_0^TX))^2/v(X)XX^T)$. 
The matrix $V^{-1}$ is called the quasi-Fisher information
attaining the minimum
 among all estimations when we use the residuals $y_i-g(\beta^Tx_i)$ to form linear combinations.
 Nevertheless,
 the above deduction is of course  under certain regularity
 conditions.

\subsection{Motivation: variance reduction by an extra ``estimation"}
 To simplify the motivation for our estimation procedures in the next subsection, we   give, for the time being, four assumptions among which the first three will be removed in the next section.  First, assume $v(\cdot)=1$ for the moment. Second,  we consider for the time being the case where  that $\beta$ is of norm one: $\|\beta\|=1$. This is because the following arguments do not make sense for general $\beta$ due to identifiability problem. This is also the reason why we need two-stage estimation that will be described in the next section. Third, to reparametrize $\beta$ when its norm is equal to one, we assume that we have a parameter $\theta$  of $p-1$ dimension so that the $p\times (p-1)$ Jacobin matrix $J(\theta)={\partial \beta}/{\partial \theta}$ exists. The details will be presented in the next section as well. Forth, the dimension of $\beta$ is greater than one, otherwise, our approach can improve nothing for the classical quasi-likelihood.

From the title of this subsection, it is a natural question to ask why and how we consider introducing an extra ``estimation" when we use the quasi-likelihood. As is known, the quasi-likelihood is motivated from the likelihood
  with given density function of the error, or the weighted  least squares
  when the distribution of the error is unknown. For simplicity of
  illustration, we will  consider the ordinary least squares (OLS)  that is equivalent to the likelihood  when the error term follows the standard normal distribution $N(0, 1)$.

\subsubsection{Motivation from estimation criterion}
  The  idea of the OLS is to search for a function $\tilde g$ with
    $\beta$ so that $\tilde g(\beta^TX)$ can fit $Y$ best in the OLS
    criterion. The basic probability theory tells us that the best $\tilde g (\beta^TX)$ should be the conditional expectation of Y given $\beta^TX.$ Note that in the GLM $E(Y|\beta_0^TX)=g(\beta_0^TX)$ with respect to the joint distribution  $F(x, y)$ of $(X, Y)$.  Thus under
    the parametric framework, a natural choice is not to search for  function $\tilde g$, while to use $g$  and then to
    consider the class of $g(\beta^TX)$
    over all $\beta.$ Thus, the OLS criterion is defined as
 $M(\beta)=E(Y-g(\beta^TX))^2$ for any $\beta$.    It is
  easy to see that the true parameter $\beta_0$ is the minimizer of
  $M(\beta)$ over all $\beta$ because
  $$E(Y-g(\beta^TX))^2=E(\varepsilon)^2+E(g(\beta_0^TX)-g(\beta^TX))^2.$$
  The derivative of $M(\beta)$  with respect to $\theta$
  leads to the population version of the estimating equation $G(\beta)$:
$$E\Big (\varepsilon g'(\beta^TX)J(\theta)^TX\Big )=E\Big([g(\beta^TX)-g(\beta_0^TX)]g'(\beta^TX)J(\theta)^TX\Big ).$$
When the expectations in the both sides are replaced by their
empirical versions that are based on the sample $\{(x_1, y_1),
\cdots, (x_n, y_n)\}$, we can obtain the estimating equation that is a special case of the
quasi-likelihood. When $\beta$ is close to $\beta_0$ at certain
rate,  $\theta$ close to $\theta_0$ that is associated with $\beta_0$,
and  $g'(\beta^TX) \approx g'(\beta_0^TX)$, the application of Taylor
expansion yields that
\begin{eqnarray}\label{mot2o}
&&E\Big (\varepsilon g'(\beta^TX)J(\theta_0)^TX\Big )\nonumber\\
&\approx&E\Big((g'(\beta_0^TX))^2J(\theta_0)^TXX^TJ(\theta_0) \Big )(\theta-\theta_0).
\end{eqnarray}
The limiting variance-covariance matrix of $\theta-\theta_0$ would be $V^{-1}$ where $${V}=E\Big ( g'(\beta_0^TX)^2J(\theta_0)^TX
 X^TJ(\theta_0)\Big ).
 $$ Further, since $\beta-\beta_0\approx J(\theta_0)(\theta-\theta_0)$, we can have that the limiting variance-covariance matrix of $\beta-\beta_0$ would be $J(\theta_0)V^{-1}J(\theta_0)^T$.

This method is entirely dependent on the parametric structure the
model under study is of.  As we mentioned above, at population level, $\beta_0$ is the minimizer (solution) of $M(\beta)$ ($G(\beta)$), and $g(\beta_0^TX)$ is the conditional
expectation $E(Y|\beta_0^TX)$. However,  at sample level, the minimizer (solution) of the empirical version of $M(\beta)$ ($G(\beta)$)   $\beta$ is not exactly equal to the true
value $\beta_0$. As such,  the search
based on $M(\beta)$ is actually a directional search along the
direction with the fixed $g$ that is not necessarily equal to the
conditional expectation of $Y$ when $\beta^TX$ is given. This results in that  the
related minimum $\beta$ may not attain the minimum that can be obtained by searching it
in a larger subspace of functions \{ $\tilde g(\beta^T\cdot)$: all $\tilde g$ and $\beta$ \}. That is,
\begin{eqnarray*} E(Y-g(\beta^TX))^2&=&
E(Y-E(Y|\beta^TX))^2+E(E(Y|\beta^TX)-g(\beta^TX))^2\\
&\ge & E(Y-E(Y|\beta^TX))^2=\inf_{\tilde g}E(Y-\tilde g(\beta^TX))^2.\end{eqnarray*} This shows that the
residual $Y-g(\beta^TX)=\varepsilon-(g(\beta^TX)-g(\beta_0^TX))$ may be
of larger variability than
$Y-E(Y|\beta^TX)=\varepsilon-(E(Y|\beta^TX)-g(\beta_0^TX))$ is of
when $\beta\not =\beta_0$.  It suggests that in the estimating
equation, using $E(Y|\beta^TX)$ may lead an estimation of $\beta_0$ with less variability in terms of  $(E(Y|\beta^TX)-g(\beta_0^TX))\equiv (E(Y|\beta^TX)-E(Y|\beta_0^TX))$ than that in terms of
$(g(\beta^TX)-g(\beta_0^TX))$. The application of Taylor expansion at
$\beta_0$ gives us an expectation that $\beta-\beta_0$ obtained from
the former would have smaller variance than this difference obtained
from the latter.  Then,  we can consider the solution of
$$
E\Big ((Y-E(Y|\beta^TX))\Big )g'(\beta^TX)J(\theta)^TX=0.
$$
Equivalently, writing  $E(Y|\beta_0^TX)=E(\varepsilon|\beta_0^TX)+E(g(\beta_0^TX)|\beta_0^TX)$,
\begin{eqnarray}\label{mot1}
&&E\Big (((\varepsilon-(E(Y|\beta_0^TX)-
g(\beta_0^TX)))) g'(\beta^TX)J(\theta)^TX\Big )\nonumber\\
&=&E\Big (((\varepsilon-E(\varepsilon|\beta_0^TX))+
(g(\beta_0^TX)-E(g(\beta_0^TX)|\beta_0^TX))) g'(\beta^TX)J(\theta)^TX\Big )\nonumber\\
&=&E\Big([E(Y|\beta^TX)-E(Y|\beta_0^TX)]
g'(\beta^TX)J(\theta)^TX\Big ).
\end{eqnarray}
Here we use the derivative $g'$ of $g$ rather than the derivative $E'(Y|\beta^TX)$ of $E(Y|\beta^TX)$ only because we try to avoid too many unknowns to be estimated. Note that in this presentation,  although both of them are zero at
population level we do not  delete $E(\varepsilon|\beta_0^TX)$ and
$(g(\beta_0^TX)-E(g(\beta_0^TX)|\beta_0^TX)).$ This is  because they are not
necessarily zero at sample  level, and  will play an important role
in variance reduction. When $\beta$ is close to $\beta_0$ at certain
rate, and also $\theta$ close to $\theta_0$ that is associated with $\beta_0$, noting that conditional expectation is a self-adjoint
operator,
and  $g'(\beta^TX) \approx g'(\beta_0^TX)$, Taylor
expansion for (\ref{mot1}) yields that
\begin{eqnarray}\label{mot2}
&&E\Big (((\varepsilon-E(\varepsilon|\beta_0^TX))) g'(\beta^TX)J(\theta_0)^TX\Big )\nonumber\\
&\approx&E\Big (\varepsilon [g'(\beta_0^TX)(X- E(J(\theta_0)^TX|\beta_0^TX))]\Big )\nonumber\\
&\approx&E\Big((g'(\beta_0^TX))^2J(\theta_0)^TXX^TJ(\theta_0) \Big )(\theta-\theta_0)\nonumber\\
&&-E\Big((g(\beta_0^TX)-E (g(\beta_0^TX)|\beta_0^TX))
g'(\beta_0^TX)J(\theta_0)^TX\Big ).
\end{eqnarray}
In terms of (\ref{mot1}) and (\ref{mot2}), we can make a comparison with the classical least squares or likelihood. Asymptotically, our approach makes the equation to have  an extra term  $E\Big ((g(\beta_0^TX)-(E(Y|\beta_0^TX)
)) g'(\beta^TX)J(\theta)^TX\Big )$ at the true value of $\beta_0.$ This term does not come up in the classical least squares because by the model structure, this term equals zero at population level, and thus, in classical framework, we do not bother ourself to ``estimate" it. More importantly, from (\ref{mot2}), we can see easily that $E\Big ((g(\beta_0^TX)-(E(Y|\beta_0^TX)
)) g'(\beta_0^TX)J(\theta_0)^TX\Big )$ is non-positively correlated to $E\Big (\varepsilon) g'(\beta^TX)J(\theta_0)^TX\Big )$. The correlation is actually between $E\Big (\varepsilon g'(\beta^TX)J(\theta_0)^TX\Big )$ and $-E\Big (E(\varepsilon|\beta_0^TX) g'(\beta_0^TX)J(\theta_0)^TX\Big )$. This helps us attain smaller variance at sample level. The details are in the following.

\subsubsection{The gain from estimating $E(Y|\beta^TX)$}
In the above presentation we note that the conditional expectation
function $E (Y|\beta^TX)$ is a
nuisance parameter to be estimated. From the approximation of (\ref{mot2}), if we can
define an estimator $\hat E (Y|\beta^TX)$ of $E (Y|\beta^TX)$ that satisfies the following requirements, an asymptotically more efficient estimation of $\beta_0$
 can be expected. First, the  estimator can make the second term on the right hand side of
(\ref{mot2}) to be asymptotically negligible. Second, its derivative at
$\beta_0^TX$, $\hat E'(Y|\beta_0^TX)$ converges to $g'(\beta_0^TX)$
at certain rate.  Third, more importantly,  an estimation of $\hat E(Y|\beta^TX))$ would
create an extra ``estimation" $\hat E(Y|\beta_0^TX))$ of
$g(\beta_0^TX)$  so that the weighted sum of $g(\beta_0^Tx_i)-\hat E(Y|\beta_0^Tx_i))$'s
is {\it non-positively correlated to} the weighted sum of $\varepsilon_i$'s for the purpose of variance reduction. From the above (\ref{mot1}) and (\ref{mot2}),
we can see that we do get an extra term about $-\hat E(\varepsilon|\beta_0^TX).$
Specifically,  such an extra estimation is not negligible, that is, $\hat E(Y|\beta_0^TX)-g(\beta_0^TX)=\hat E(\varepsilon|\beta_0^TX)+(\hat E(g(\beta_0^TX)|\beta_0^TX)-g(\beta_0^TX))$ is not negligible. Technically, we will prove that $(\hat E(g(\beta_0^TX)|\beta_0^TX)-g(\beta_0^TX))$ is negligible. Thus, the gain we can have is from $\hat E(\varepsilon|\beta_0^TX)$. Without confusion, here $\hat E$ means estimation  of conditional mean and  of unconditional mean for its appearance in different places.  Note that $\hat E(\varepsilon g'(\beta_0^TX)J(\theta_0)^TX)$ is  asymptotically normal and helps us achieve the Fisher information. When  $\hat E\Big ((g(\beta_0^TX)-\hat E(Y|\beta_0^TX))g'(\beta_0^TX)J(\theta_0)^TX\Big )\approx  - \hat E\Big ( \hat E(\varepsilon|\beta_0^TX)g'(\beta_0^TX)J(\theta_0)^TX\Big )  $ is {\it non-positively related to} $\hat E(\varepsilon g'(\beta_0^TX)J(\theta_0)^TX)$,   we will have chance to have a smaller variance than the Fisher information. This is just the case with our approach. The covariance between these two terms is  asymptotically equal to $$-Cov:=-E(g'(\beta_0^TX)^2E(J(\theta_0)^TX|\beta_0^TX)E(J(\theta_0)X^T|\beta_0^TX)).$$ Regardless of some technical details for nonparametric estimation, this calculation can be easily justified  from (\ref{mot1}) and (\ref{mot2}).  Interestingly,  We can also derive that $\hat E\Big ((g(\beta_0^TX)-\hat E(Y|\beta_0^TX))\Big )g'(\beta_0^TX)J(\theta_0)^TX$ has a limiting variance-covariance matrix  equal to $Cov$. In other words, the joint distribution of
$$(\hat E(\varepsilon g'(\beta_0^TX)J(\theta_0)^TX), \, \,  \hat E\Big ((g(\beta_0^TX)-\hat E(Y|\beta_0^TX))g'(\beta_0^TX)J(\theta_0)^TX)\Big )
$$
is asymptotically normal with the covariance matrix
$$\left (
\begin{array}{cc}
E(g'(\beta_0^TX)^2J(\theta_0)^TXX^TJ(\theta_0)) & -Cov\\
-Cov&
Cov\\
\end{array}\right ).$$
Thus, the sum of these two terms is asymptotically normal with the variance-covariance matrix
$$\tilde Q=V-Cov.$$
We can also obtain it from (\ref{mot2}). Therefore, an extra estimation  $\hat E(Y|\beta^TX)$ of
 $E(Y|\beta^TX)$ automatically introduces an extra ``estimation"  $\hat E(Y|\beta_0^TX)$
 of $E(Y|\beta_0^TX)=g(\beta_0^TX)$ that can help us achieve the reduction
 of  variance.

By comparing the two sides of (\ref{mot2}), the limiting variance-covariance
 matrix of $\theta-\theta_0$ would be $V^{-1}{\tilde Q}V^{-1}$ where
 $${\tilde Q}=E\Big ( g'(\beta_0^TX)^2J(\theta_0)^T(X-
 E(X|\beta_0^TX))(X^T- E(X^T|\beta_0^TX))J(\theta_0)\Big ),
 $$ and recalling its definition
 $${V}=E\Big ( g'(\beta_0^TX)^2J(\theta_0)^TX
 X^TJ(\theta_0)\Big ).
 $$
 It is clear that $V-{\tilde Q}=Cov$ is a non-negative semidefinite
 matrix
 and then $V^{-1}{\tilde Q}V^{-1}$ is smaller than or equal to  $V^{-1}.$ From this, and $\beta-\beta_0\approx J(\theta_0)(\theta-\theta_0)$, we may obtain that the limiting variance-covariance matrix would be $J(\theta_0)V^{-1}{\tilde Q}V^{-1}J(\theta_0)^T$. This matrix is smaller than or equal to $J(\theta_0)V^{-1}J(\theta_0)^T.$

The above "super-efficiency" can also be found to have a root from estimation problems with nuisance parameters as we mentioned in Section~1. For instance, suppose that we have a model with two sets of parameters: $\beta_0$, the parameters of interest; and $\lambda$, nuisance parameters. Piece (1982) found that an estimation with estimated $\lambda$ can sometimes be more efficient than an ``estimation" with the true value of $\lambda$. In our setting, we regard the conditional expectation $E(Y|\cdot)$ as a nuisance parameter that is actually a pseudo-nuisance parameter for the models, and the link function  $g(\cdot)$ is regarded as its true value. In the classical methodologies,  the true value is used in estimations, whereas in our method,  an estimation of $E(Y|\cdot)$ is plugged in.

\section{Two-stage estimation procedures}\label{two-pahse}
From the above explanations about the least squares and likelihood, we can see the intrinsic difference of our approach from the classical estimations.
This motivates us  to use the  following quasi-likelihood for the
purpose of variance reduction:
\begin{eqnarray}\label{gee4}
SG(\beta)=\sum_{i=1}^n(y_i-\hat E(Y|\beta^Tx_i))
g'(\beta^Tx_i)x_i/v(x_i)I_n(x_i),
 \end{eqnarray}
 where $\hat E(Y|\beta^TX)$ is an estimator of $E(Y|\beta^TX)$,
  $I_n(x_i) = I
  \{\hat f_{\hat \beta}(\hat \beta^Tx_i)>c_0\}$ for some positive value $c_0$, and $\hat f$ and $\hat \gamma_0$ are respectively estimators of the density function of $\beta_0^TX$ and $\beta_0$ in (\ref{tranca}) of Section~3 where $\beta_0$ is replaced by $\gamma_0$, the direction of $\beta_0$. We can see the details in Section~3 below. The truncation $I_n(x_i)$
  is employed here for technical
  purpose to handle the boundary points particularly when we use local smoothing method to estimate the nonparametric
  regression function $E(Y|\beta^TX)$, see Xia, and H\"{a}rdle(2006).

{\it However,  the above arguments are only feasible in the case where the forth assumptions stated in the beginning of this subsection satisfy}. The constraint $\|\beta_0\|=1$ is intrinsic for our approach because for general $\beta$, $E(Y|\beta^TX)=E(Y|c\beta^TX)$ for any constant $c$ and then  the uniqueness of  solution $\beta$ cannot be
guaranteed. To deal with  this
identifiability problem, we in the following propose two-stage
estimation procedures.

\subsection{Re-parametrization of the parameter}
We note that although $\beta_0$ is not identifiable in the above procedure, its direction
$\beta_0/\|\beta_0\|$ can be so when we apply the idea described in the previous subsection. As such, we can either estimate the
direction first or its scale $\|\beta_0\|$ first. Write
$\alpha_0=\|\beta_0\|$  and $\gamma_0=\beta_0/\|\beta_0\|$. It is
easy to see that $\gamma_0$ is identifiable when we rewrite the
original model as $Y=E(Y|\gamma_0^TX)+\varepsilon$. We can work on
estimating $\gamma_0$ by using a similar estimating equation of
(\ref{gee4}) and $\alpha_0$ can be estimated by the classical
quasi-likelihood. However, we note that $\gamma_0$ is on the
boundary of the unit sphere surface and then we cannot directly
derive the scores in the quasi-likelihood. To deal with this, a
re-parametrization is necessary. A popular re-parametrization is the
"remove-one-component" method on $\beta$ as in Yu and Ruppert(2002)
and then Wang, Xue, Zhu and Chong (2010) and Chang, Xue and Zhu
(2010).
 Without loss of generality, we may assume that the true
parameter $\gamma_0$ has a positive component, say $\gamma_{0r}>0$
for $\gamma_{0} = (\gamma_{01},\ldots,\gamma_{0p})^T$ and $1\leq
i\leq n$. For $\gamma = (\gamma_{1},\ldots,\gamma_{p})^T$, let
$\gamma^{(r)}=(\gamma_{1},\ldots,\gamma_{r-1},\gamma_{r+1},\ldots,\gamma_{p})^T$
be the $p-1$ dimensional parameter vector after removing the $r$th
component $\gamma_r$ in $\gamma$. We may write
\begin{eqnarray}\label{betaexp}
\gamma =
\gamma(\gamma^{(r)})=(\gamma_{1},\ldots,\gamma_{r-1},(1-\parallel\gamma^{(r)}\parallel^2)^{1/2},\gamma_{r+1},\ldots,\gamma_{p})^T.
\end{eqnarray}
The true parameter $\gamma^{(r)}_0$ must satisfy the constraint
$\parallel\gamma^{(r)}_0\parallel^2<1$. Thus, $\gamma$ is infinitely
differential in a neighborhood of $\gamma^{(r)}_0$, and the Jacobian
matrix is
\begin{eqnarray}\label{jacobin}
J(\gamma^{(r)})=\frac{\partial\gamma}{\partial\gamma^{(r)}}=(\delta_1,\ldots,\delta_p)^T,
\end{eqnarray}
where $\delta_s(1\leq s\leq p)$ satisfy that $\delta_s=e_s$ for
$1\leq s<
r$,$\delta_r=-(1-\parallel\gamma^{(r)}\parallel^2)^{-1/2}\gamma^{(r)}$,$\delta_s=e_{s-1}$
for $r+1\leq s\leq p$. Here $e_s$ is a $p-1$ dimensional unit vector
with $s$th component $1$. For this re-parametrization, we can also
prove that $J(\gamma^{(r)})$ is orthogonal to $\gamma$. However, the
columns within $J(\gamma^{(r)})$ are not orthogonal. Let
$A(\gamma^{(r)})=(\gamma, J(\gamma^{(r)}))$. We can have
$$A(\gamma^{(r)})^TA(\gamma^{(r)})
=\left(
        \begin{array}{cc}
          1 & {\bf 0}^T \\
          {\bf 0} & J(\gamma^{(r)})^TJ(\gamma^{(r)}) \\
        \end{array}
      \right)$$ although it  is not a diagonal matrix.
        Also, we can easily prove that  $A(\gamma^{(r)})$
         is invertible, and its inverse is
\begin{eqnarray}A(\gamma^{(r)})^T&=&
\left(
        \begin{array}{cc}
          1 & {\bf 0}^T \\
          {\bf 0} & J(\gamma^{(r)})^TJ(\gamma^{(r)}) \\
        \end{array}
      \right)A(\gamma^{(r)})^{-1}, \nonumber\\
       A(\gamma^{(r)})&=& (A(\gamma^{(r)})^T)^{-1}
\left(
        \begin{array}{cc}
          1 & {\bf 0}^T \\
          {\bf 0} & J(\gamma^{(r)})^TJ(\gamma^{(r)}) \\
        \end{array}
      \right).
  \end{eqnarray}

  An alternative is re-parametrization by the
 polar coordinate system. to simplify the presentation, we assume with no loss of generality that the $p$-th component of $\gamma$ is positive. That is,
letting $\theta=(\theta_1, \cdots, \theta_{p-1})^T$ whose domain is
the $(p-1)$-dimensional subspace $(0, \pi)^{p-2}\times (0, \pi/2)$,
\begin{eqnarray}\label{gamma1}\gamma=\gamma(\theta)=
\Big (\prod_{i=1}^{p-1}\cos(\theta_i), \sin
(\theta_1)\prod_{i=2}^{p-1}\cos(\theta_i),
\sin(\theta_2)\prod_{i=3}^{p-1}\cos(\theta_i) \cdots,
\sin(\theta_{p-1}) \Big )^T,
  \end{eqnarray}
where $\prod_{i=p}^{p-1}\cos(\theta_i) $ is defined as $1$. Otherwise, $\sin(\theta_{p-1}$ will be  the $r$-th component that is positive.
 The Jacobian matrix $J(\theta)=\partial \gamma(\theta)/\partial \theta$
 is of a special structure: the $i$-th column is the derivative of
 $\gamma(\theta)$ about $\theta_i$, in the first $i$ components,
 $\sin(\theta_i)$  and $\cos (\theta_i)$ are respectively replaced
 by $\cos(\theta_i)$ and $-\sin(\theta_i)$, and  the other components
 in the column are equal to zero because the corresponding components do not have $\theta_i$.
 It is easy to prove that all the columns are orthogonal to each other.
  Furthermore, $J(\theta)$ is orthogonal to $\gamma(\theta)$.  A brief
  justification can be done by induction. When $p=2$, the conclusion is
  clearly true. Assume that when $p=m-1$ the conclusion is true.
  When $p=m$,
   $\gamma_m=
    (\cos(\theta_m)\gamma_{m-1}^T, \sin(\theta_m))^T$.
    For the derivative about $\theta_i$ with
    $i\le m-1$,
 \begin{eqnarray}\label{jacob1}
 \partial \gamma_m/\partial \theta_i&=&
 (\cos(\theta_m)\partial \gamma_{m-1}^T/\partial
 \theta_i , 0)^T\nonumber
 \end{eqnarray}
 Then for any $1\le j<i$,
 $$(\partial \gamma_{m}^T/\partial \theta_i)(\partial
 \gamma_{m}/\partial \theta_j)=(\partial \gamma_{m-1}^T
 /\partial \theta_i)(\partial \gamma_{m-1}/\partial \theta_j)
 \cos(\theta_m)^2=0.$$
  When $i=m$, $\partial \gamma_{m}/\partial \theta_m=
(\sin(\theta_m) \gamma_{m-1}^T, \cos(\theta_m))^T,$
   we have for any $j<m$
\begin{eqnarray}\label{jacob2}
  &&(\partial \gamma_{m}^T/\partial \theta_j)( \partial \gamma_{m}/\partial \theta_m )\nonumber\\
  &=& (\sin(\theta_m)\partial \gamma_{m-1}^T/\partial \theta_j, 0)
  (\cos(\theta_m) \gamma_{m-1}^T, -\sin(\theta_m))^T\nonumber\\
  &=&  \cos(\theta_m)\sin(\theta_m)(\partial \gamma_{m-1}^T/\partial \theta_j) \gamma_{m-1}
=0.\nonumber
\end{eqnarray}
But  $A(\theta):=(\gamma(\theta), J(\theta))$  is not an orthogonal
matrix although
$$A(\theta)^TA(\theta)
=\left(
        \begin{array}{cc}
          1 & {\bf 0}^T \\
          {\bf 0} & J(\theta)^TJ(\theta) \\
        \end{array}
      \right).$$  is a diagonal matrix as not all elements are equal to $1$.
        Also, we can easily prove that
\begin{eqnarray}A(\theta)^T=
\left(
        \begin{array}{cc}
          1 & {\bf 0}^T \\
          {\bf 0} & J(\theta)^TJ(\theta) \\
        \end{array}
      \right)A(\theta)^{+}, \,  \, A(\theta)= (A(\theta)^T)^{+}
\left(
        \begin{array}{cc}
          1 & {\bf 0}^T \\
          {\bf 0} & J(\theta)^TJ(\theta) \\
        \end{array}
      \right).
  \end{eqnarray}
where $A^{+}$ is the Moore-Penrose inverse of matrix $A$ in case $A$
is not invertible. This will be useful in the later analysis.

We can see that the two re-parametrizations have some similar
properties. But, the former is easier to compute and
$A(\gamma^{(r)})$ is always invertible whereas the latter has nicer
orthogonality
structure.   Throughout the rest of the present paper, 
 we assume with no loss of generality that $A(\theta)$ is invertible in
the following, $\theta$ can be either $\gamma^{(r)}$ in the former
or $\theta$ in the latter.

\subsection{Estimation procedures}
Now we are in the position to present
estimation procedures. Note that $E(Y|\gamma_0^TX)
=g(\alpha_0\gamma_0^TX)$. The derivative of $g(\alpha
\gamma(\theta)^TX)$ about $\theta$ is $g'(\alpha
\gamma(\theta)^TX)\alpha J(\theta)^T X$. Further, as the estimation
will involve nonparametric smoothing,  the boundary effect needs to
be dealt with, we will in the estimating equation trim off some
boundary points in the following procedures.

 \noindent{\bf Procedure 1}. \\
{\it Step 1.} Obtain an initial estimator of $\beta_0$, $\hat
\beta_I$ by the estimating equation of (\ref{gee1}),
and then let $\hat \gamma_I=\hat \beta_I/\|\hat \beta_I\|$. Check the values of all components of $\hat \gamma_I$. Select an $r$ with $1\le r\le p$ such that $|\hat \gamma_I^{(r)}|=\max_{1\le l\le p}|\hat \gamma_I^{(l)}|.$  As $\hat \gamma_I^{(r)}$ is not necessary to be positive, we then define $\hat \alpha_I=sign(\hat \gamma_I^{(r)})\|\hat \beta_I\|.$  As such, we can define the domain $\Gamma(\Theta)$ of $\gamma(\theta)$ with the constraint that the $r$th component of $\gamma(\theta)$ is positive.\\
{\it Step 2.} Estimate $\theta_0$ by the solution $\hat \theta$ of
\begin{eqnarray}\label{ngee1}
SG(\theta)=\sum_{i=1}^n(y_i-\hat E(Y|\gamma(\theta)^Tx_i)) g'(\hat
\alpha\gamma(\theta)^Tx_i)J(\theta)^Tx_i/v(x_i)\tilde{I}_n(x_i)=0.
 \end{eqnarray}
over the domain $\Gamma(\Theta)$, where $\tilde{I}_n(x_i)$ is a truncation function in (\ref{gee4}).\\ 
 {\it Step 3.} Obtain the final estimator of  $\beta$ as $\hat \beta_F=
 \hat \alpha_I  \gamma(\hat \theta)$.

 Another procedure uses an iterative algorithm with one more new Step benefiting from the variance reduction estimation of $\gamma$ above.

 \noindent{\bf Procedure 2}.
 \\
{\it Step 1.} Obtain an initial estimator of $\beta_0$, $\hat
\beta_I$ by the estimating equation of (\ref{gee1}),
and then let $\hat \gamma_I=\hat \beta_I/\|\hat \beta_I\|$. Define $\hat \alpha_I=\|\hat \beta_I\|.$  \\
{\it Step 2.} Estimate $\theta_0$ by the solution $\hat \theta$ of
\begin{eqnarray}\label{ngee1}
SG(\theta)=\sum_{i=1}^n(y_i-\hat E(Y|\gamma(\theta)^Tx_i)) g'(\hat
\alpha\gamma(\theta)^Tx_i)J(\theta)^Tx_i/v(x_i)\tilde{I}_n(x_i)=0.
 \end{eqnarray}
over all $\theta)$, where $\tilde{I}_n(x_i)$ is a truncation function in (\ref{gee4}).\\ 
{\it Step 3.} Estimate $\alpha_0$ by the solution $\hat \alpha_F$ of
\begin{eqnarray}\label{ngee2}
G(\alpha)=\sum_{i=1}^n(y_i-g(\alpha \gamma(\hat \theta)^Tx_i)) g'(
\alpha \gamma(\hat \theta)^Tx_i) \gamma(\hat
\theta)^Tx_i/v(x_i)\tilde{I}_n(x_i)=0.
 \end{eqnarray}\\
 {\it Step 4.} Obtain the final estimator of  $\beta$ as $\hat \beta_F= \hat \alpha_F \gamma(\hat \theta)$.


\begin{rema} The two procedures are different mainly for Step~1 and Step~3. When we assume that a component of $\gamma$ is positive, the scalar $\alpha$ is then not necessary to be positive. To avoid this identification issue about sign, we define an $\hat \alpha_I$ with $sign$ function in Step~1 of Procedure~1. But this is no need for Procedure~2 as after having an estimation of $\gamma$, we re-estimate the scalar $\alpha$, with which, we can adoptively have the sign of $\alpha$ in the final estimation.
\end{rema}

\begin{rema} Trimming boundary points off   is a typical technique in nonparametric estimation when local smoothing is applied.
We here simply use $\tilde{I}_n(x_i)$. Similar trimming can be found
in Xia and H\"ardle (2006), and Xia, H\"ardle and Linton (2009) who
used a smooth function instead of the indicator function. To achieve
the results without such a function appeared in the limiting
variance, that is, $\tilde{I}_n(x_i)$ goes to $1$, we should let
$c_0$ tend to zero as $n$ goes to infinity. It is achievable as long
as $c_0$ tends to zero at certain rate not too fast, see Xia and
H\"ardle (2006), and Xia, H\"ardle and Linton (2009) for more
details. Thus, we will not discuss the choice of $c_0$ in detail.

\end{rema}

\subsection{Plug-in Estimation of $E(Y|\gamma^TX)$}

In the estimation procedures above, we need a plug-in estimation for
$E(Y|\gamma^TX)$. Here we consider local linear smoother. Suppose
that $(x_i,y_i)$,$1\leq i\leq n$ are independent identically
distributed (iid) from model (\ref{glm2}).  To estimate
nonparametrically $E(Y|\gamma^TX)$, the local linear smoother with
the linear weighted functions $\{W_{ni}(z,\gamma,h);1\leq i\leq n\}$
is defined in the following. The local linear estimators of
$E(Y|\gamma^TX)$ and its derivative are the pair $(a,b)$ minimizing
the weighted sum of squares (see Fan and Gijbels 1996)
\begin{eqnarray}
\sum^n_{i=1}\left[y_i-a-b(x_i^T\gamma-z)\right]^2K_h(x_i^T\gamma-z),
\end{eqnarray}
where $K_h(\cdot)=K(\cdot/h)/h$ with $K$ being a symmetric kernel
function on the real line and $h=h_n$ being a bandwidth.

Without confusion, we use the notation $\hat{g}(z;\gamma,h):=\hat
E(Y|\gamma^TX=z)$. $\hat{g}(z;\gamma,h)=\hat a$ is the resulting
estimator. Via a simple calculation, we have
\begin{eqnarray}
\displaystyle\hat{g}(z;\gamma,h)=\sum^n_{i=1}W_{ni}(z;\gamma,h)y_i,
\end{eqnarray}
where, for $1\leq i\leq n$,
\begin{eqnarray}
W_{ni}(z;\gamma,h)=\frac{U_{ni}(z;\gamma,h)}
{\sum^n_{i=1}U_{ni}(z;\gamma,h)}
\end{eqnarray}
with
\begin{eqnarray}
U_{ni}(z;\gamma,h)&=&[S_{n,2}(z;\gamma,h)-(x_i^T\gamma-z)S_{n,1}(z;\gamma,h)]
K_h(x_i^T\gamma-z)\nonumber,\\
S_{n,r}(z;\gamma,h)&=&\frac{1}{n}\sum^n_{i=1}(x_i^T\gamma-z)^rK_h(x_i^T
\gamma-z),\quad r=0,1,2.\nonumber
\end{eqnarray}
For the density function  ${f}_{{\gamma}}({\gamma}^TX)$, we can use the following estimator \begin{eqnarray}\label{tranca}\hat{f}_{\hat{\gamma}}(\hat{\gamma}^Tx_i)=
 (nh^2)^{-1}{\sum^n_{i=1}U_{ni}(\hat{\gamma}^Tx_i;\hat \gamma,h)},\end{eqnarray} and  $\hat \gamma=\hat \beta_I/\|\hat \beta_I \|$.
\section{Main results}\label{sec4}
 By the classical quasi-likelihood for the GLM, we can have an initial
  estimator $\hat \beta_I$ of $\beta_0$ that is root-$n$ consistent.
 It implies that $\hat \gamma_I$ is also root-$n$ consistent to $\gamma_0$.
  With this
  pilot estimation, we can then consider a
neighborhood of $\hat \gamma_I$  which contains $\gamma_0$ to define
our resulting estimator. Let
$\mathcal{B}_n=\{\gamma:\parallel\gamma-\hat
\gamma_I\parallel=Bn^{-1/2}\}$ where $B$ is some positive constant.

\begin{theo}\label{theo1}
Suppose that conditions {\rm C1--C6}  hold. Let the solution of
(\ref{ngee1}) be $\gamma(\hat \theta)$. Assume further that $f(x)>c_0$ and
$f_{\gamma}(\gamma^Tx)>c_0$ for all $x$ in its support. We have
\begin{eqnarray}
\sqrt{n}(\gamma(\hat
\theta)-\gamma_0)\stackrel{\mathcal{L}}{\longrightarrow}N(0,
J(\theta_0){\tilde V}^{-1}Q_1{\tilde
V}^{-1}J(\theta_0)^T/\alpha_0^2), \nonumber
\end{eqnarray}
where $Q_1=:J(\theta_0)^T{\tilde
Q}J(\theta_0)=E\{|\varepsilon|^2g'(\beta_0^TX)^2J(\theta_0)^T
[X/v(X)-E(X/v(X)|\gamma_0^TX)][X/v(X)-E(X/v(X)|\gamma_0^TX)]^T
J(\theta_0)\}$, $\tilde{V}= E\{[g'(\beta_0^TX)]^2 J(\theta_0)^T
XX^T/v(X) J(\theta_0)\} $. 
\end{theo}
\begin{rema}\label{rema1.1} The above result is obtained under the condition $f(x)>c_0$ and
$f_{\gamma}(\gamma^Tx)>c_0$ for all $x$ in its support for
simplicity. When these conditions do not hold, in (\ref{ngee1}),
$\tilde I_n(x_i) = I\{
 \hat{f}_{\hat{\gamma}}(\hat{\gamma}^Tx_i)>c_0\}$ is not equal to $1$, where, recalling its definition in (\ref{tranca}) of the previous subsection,  $\hat{f}_{\hat{\gamma}}(\hat{\gamma}^Tx_i)=
 (nh^2)^{-1}{\sum^n_{i=1}U_{ni}(\hat{\gamma}^Tx_i;\hat \gamma,h)}$, $\hat \gamma=\hat \beta_I/\|\hat \beta_I \|$. From the proof of Lemma~\ref{lemma2}, we can see that  $ \hat{f}_{\hat{\gamma}}(z)$ converges to a function $\mu_2f^2_{\gamma_0}(z)$ where
 $\mu_2$  to be specified in Appendix is related to the kernel function in local linear smoother. In the place of $X/v(X)$ we should use $X/v(X)\tilde I(X)$ where $\tilde I(X)$ is the limit of $\tilde I_n(X)$ as $n$ tends to infinity. The above variance in the theorem can be regarded as the limit  as $c_0$  go to zero at a certain rate.
  All the
results below will be with  $\tilde I_n(X)=1$ for simplicity.
Otherwise, the proof will be similar with more tedious computation,
and the result will  be substituted by
$Q_1=E\{|\varepsilon|^2g'(\beta_0^TX)^2J(\theta_0)^T [X{\tilde
I}(X)/v(X)-E(X{\tilde I}(X)/v(X)|\gamma_0^TX)][X{\tilde
I}(X)/v(X)-E(X{\tilde I}(X)/v(X)|\gamma_0^TX)]^T
J(\theta_0)\}=:J(\theta_0)^T{\tilde Q}J(\theta_0)$, and \newline
$\tilde{V}= E\{[g'(\beta_0^TX)]^2 J(\theta_0)^T
XX^T {\tilde I}(X)/v(X)J(\theta_0)\} $ accordingly throughout this whole paper. 
\end{rema}

\begin{rema}\label{rema1.2}
For the single-index model, Wang, Xue, Zhu and Chong (2010)
 and Chang, Xue and Zhu (2010) also used estimating equations.
 However, their approaches are based on the least squares and the
 estimating equations are the derivatives of the least squares
  criterion. Thus, although their estimation is asymptotically
  more efficient than existing ones such as
  H\"ardle,  Hall and Ichimura (1993) and Xia and H\"ardle (2006),
  it cannot have optimal properties what the quasi-likelihood can have.
  More importantly, when $v(X)$ is of a structure $v(\beta_0^TX)$ for
  which we will have discussions in Section~4, we cannot usually define
  consistent estimation from
  the estimating equations derived from
  weighted least
  squares, see Heyde(1997, page 4). Also, as in our setting, we can directly use the derivative
  $g'(\beta^TX)$   as $g'$ is known. This is also different from the single-index model with unknown $g'$. The estimation procedure then has less nonparametric smoothing involved.
\end{rema}
To make a comparison with the one derived by the classical
quasi-likelihood, we give a proposition.
\begin{prop} \label{prop}
Let $\hat \beta$ is the quasi-likelihood  estimator of $\beta_0$,
and let $\hat \gamma=\hat \beta/\|\hat \beta\|$, $\hat \alpha=\|\hat
\beta\|$. Then the limiting variance of $\sqrt n (\hat
\alpha-\alpha_0)$ is $\gamma_0^TV^{-1}\gamma_0$ and the limiting
variance of $\sqrt n (\hat \beta/\|\hat \beta\|-\gamma_0)$ is $({\bf
0}_p, J(\theta_0))(A(\theta_0)^TVA(\theta_0))^{-1}({\bf 0}_p,
J(\theta_0)^T/\alpha_0^2$. This variance is greater than or equal to $
J(\theta_0){\tilde V}^{-1}Q_1{\tilde
V}^{-1}J(\theta_0)^T/\alpha_0^2$ in Theorem~\ref{theo1}. Here $V =
E[g'^2(\beta_0^TX)XX^T/v(X)]$ and ${\bf 0}_p$ is a zero vector of $p$ dimension. 
\end{prop}

\begin{rema}\label{rema1.3}
This proposition also shows that when the norm of $\beta_0$ is one, that is, $\beta_0\equiv \gamma_0$, the limiting variance covariance matrix is smaller or equal to the quasi-Fisher information (or the Fisher information in the exponential distribution case).
In general cases, this proposition helps us obtain the final estimator $\hat \alpha
\gamma(\hat \theta)$ with smaller limiting variance as follows.
\end{rema}


\begin{theo} \label{theo2}
Under the conditions of Theorem~\ref{theo1}, when {\bf Procedure} 1 is
applied, we have
\begin{eqnarray}
\sqrt n (\hat \alpha \gamma(\hat \theta)-\beta_0)
\stackrel{\mathcal{L}}{\longrightarrow}N(0,
A(\theta_0)WA(\theta_0)^T),
\end{eqnarray} where $W=\left(
        \begin{array}{cc}
          \gamma_0^TV^{-1}\gamma_0 & {\bf 0}_{p-1}^T \\
          {\bf 0}_{p-1} & \tilde{V}^{-1}Q_1\tilde{V}^{-1}\\
        \end{array}
      \right),$
$\tilde{V}$ and $V$ are respectively defined in Theorem~\ref{theo1} and
Proposition~\ref{prop}.
\end{theo}

For the estimator obtained by {\bf Procedure}~2, we state the following theorem.
\begin{theo} \label{theo3}
Under the conditions of Theorem~\ref{theo1}, when {\bf Procedure}~2 is
adopted, we have
\begin{eqnarray}
\sqrt n (\hat \alpha_F \gamma(\hat \theta)-\beta_0)
\stackrel{\mathcal{L}}{\longrightarrow}N(0,
A(\theta_0)W_1A(\theta_0)^T),
\end{eqnarray} where $W_1=\left(
        \begin{array}{cc}
          (\gamma_0^TV\gamma_0)^{-1}+
          \frac {\gamma_0^TVJ(\theta)\tilde{V}^{-1}Q_1\tilde{V}^{-1}J(\theta_0)^TV\gamma_0}
          {(\gamma_0^TV\gamma_0)^2}& {\bf 0}_{p-1}^T \\
          {\bf 0}_{p-1} & \tilde{V}^{-1}Q_1\tilde{V}^{-1} \\
        \end{array}
      \right),$ and $W\ge W_1$.
\end{theo}

The proof in the appendix shows that the first element of $W$ is
greater than or equal to the corresponding element of  $W_1$. Then
$W\ge W_1$. This indicates that when we use the estimator $\hat
\gamma$ as a plug-in for estimating $\alpha_0$, the resulting
estimation does benefit from its smaller variance. This suggests
that the estimation of $\alpha_0$  directly from the classical
quasi-likelihood is, at least asymptotically, also worse than the
new method. We are now in the position to compare the matrices
involved.

\begin{theo} \label{theo4}
Let $V_2=A(\theta_0)^TVA(\theta_0)$. Then
$V^{-1}=A(\theta_0)V_2^{-1}A(\theta_0)^{T}.$ We have\\
\quad 1. All the elements on the diagonal of $W$ (or $W_1$) are
 smaller than or equal to the corresponding elements of $V^{-1}.$\\
\quad 2.
\begin{eqnarray*}
trace(V^{-1})&\ge & trace ( A(\theta_0)WA(\theta_0)^T )(or \, \,
trace(A(\theta_0)W_1A(\theta_0)^T).
\end{eqnarray*}
%
\end{theo}

\begin{rema}

In the proof in the Appendix, we can see that when
$J(\theta_0)^TX/v(X)$ is uncorrelated to
$(g'(\beta_0^TX))^2\gamma_0^TX$, $V_2^{-1}\ge W ( = W_1)$. Further,
when a stronger assumption holds:
$E(J(\theta_0)^TX/v(X)|\gamma_0^TX)=0$, we have  $Q_1=\tilde{V}$.
Thus, $V_2^{-1}= W (= W_1),$ and
$V^{-1}=A(\theta_0)WA(\theta_0)^{T}$. This shows that the classical
quasi-likelihood can be asymptotically equally efficient to the
two-stage quasi-likelihood when $J(\theta_0)^TX$ is uncorrelated to
$\gamma_0^TX$ in certain sense. This is the case when $X$ follows a
spherically
symmetric distribution. Otherwise, the variance obtained by our method
is smaller than that by the classical quasi-likelihood in the sense stated in the theorems.
From these phenomenon, we can see that
the estimation efficiency of the new method may benefit from
the correlations among the predictors in certain sense. 

%
\end{rema}

\section{Further discussions}\label{sec3}

As this is the first step towards a general idea  of using semiparametric approaches to study parametric problems for the GLM, there are many potential issues worthy of  further exploration.

1. In the previous sections, the conditional variance of $\varepsilon$
given $X$ is assumed to be known. In many GLM models, it is of a
structure $v(\beta_0^TX)$ where $v(\cdot)$ is given up to a dispersion constant and $\beta_0$ is the
parameter of interest in the regression function. In this case, when
we use the quasi-likelihood to search for a solution, this unknown
$\beta_0$ should be either replaced by an initial estimator or be estimated by
the same quasi-likelihood. The estimating equation for the initial
estimation $\hat \beta$ is
\begin{eqnarray}\label{gee3}
G(\beta)=\sum_{i=1}^n(y_i-g(\beta^Tx_i))g'(\beta^Tx_i)x_i/v(\beta^Tx_i)=0,
 \end{eqnarray}
assuming that $v$ is given, and in {\bf Procedure}~1, we change  the estimating equation to obtain the solution of $\theta_0$ by
 \begin{eqnarray}\label{ngee3}
SG(\gamma(\theta))=\sum_{i=1}^n(y_i-\hat E(Y|\gamma(\theta)^Tx_i))
g'(\hat \alpha\gamma(\theta)^Tx_i)J(\theta)^Tx_i/v(\hat
\alpha\gamma(\theta)^Tx_i)\tilde{I}_n(x_i)=0,\nonumber\\
 \end{eqnarray}\\
 or simply use $\hat \beta$ in $v(\cdot).$
From the motivation in Section~2 and the proof of
Theorem~\ref{theo1}, it is easy to see that the estimator of
$\beta_0$ that is obtained by the modified {\bf Procedure}~1 with (\ref{ngee3}) is asymptotically not affected by any consistent
estimation $\hat \alpha\gamma(\hat \theta)$ ( or $\hat \beta$) in
$v(\hat \alpha\gamma(\hat \theta)^T\cdot)$ (or $v(\hat \beta^T\cdot)$), and similarly for
{\bf Procedure}~2 when we use
 \begin{eqnarray}\label{ngee4}
G(\alpha)=\sum_{i=1}^n(y_i-g(\alpha \gamma(\hat \theta)^Tx_i)) g'(
\alpha \gamma(\hat \theta)^Tx_i) \gamma(\hat \theta)^Tx_i/v(\alpha
\gamma(\hat \theta)^Tx_i)=0.
 \end{eqnarray}
 Then  the asymptotic properties are very similar  to those in the previous theorems as if $v(\cdot)$ were given.
Actually, the above can be more general with a general function $w(\beta^TX)$ in the lieu of $ g'(
 \beta^TX)/v( \beta^TX),$ and the technical proof can be very similar.

Further, when $v(X)$ is of a single-index structure as
$v(\alpha_0\gamma_0^TX)$ with unknown function $v(\cdot)$,  we can apply
nonparametric smoothing to the squared residuals $(y_i-g(\hat
\beta^Tx_i))^2$ to construct, for any fixed $\gamma$, an estimation
$\hat v(\hat \alpha\gamma^TX)$. This is possible because this function is the
conditional expectation of $\varepsilon^2$ given $\gamma^TX$ where
$\hat \beta$ is an initial estimator obtained by the classical
quasi-likelihood  as
\begin{eqnarray}\label{gee5}
G(\beta)=\sum_{i=1}^n(y_i-g(\beta^Tx_i))g'(\beta^Tx_i)x_i=0,
 \end{eqnarray}
As is well known, this $\hat \beta$ is not asymptotically efficient,
but still asymptotically normal and $\hat v(\cdot)$ is of a nonparametric
convergence rate (see, Fan and Gijbels 1996). This helps to define our final
estimation of $\beta_0$. In {\bf Procedure}~1, we change  the estimating
equation to obtain the solution of $\theta_0$ by
 \begin{eqnarray}\label{ngee5}
SG(\gamma(\theta))=\sum_{i=1}^n(y_i-\hat E(Y|\gamma(\theta)^Tx_i))
g'(\hat \alpha\gamma(\theta)^Tx_i)J(\theta)^Tx_i/{\hat v}(\hat
\alpha\gamma(\theta)^Tx_i)\tilde{I}_n(x_i)=0\nonumber\\
 \end{eqnarray}
 or in (\ref{ngee5}) we simply use $\hat \beta$ in lieu of  $\hat
\alpha\gamma(\theta)$ in the estimated variance function ${\hat v}(\hat
\alpha\gamma(\theta)^T\cdot).$
In {\bf Procedure}~2, we use
 \begin{eqnarray}\label{ngee6}
G(\alpha)=\sum_{i=1}^n(y_i-g(\alpha \gamma(\hat \theta)^Tx_i)) g'(
\alpha \gamma(\hat \theta)^Tx_i) \gamma(\hat \theta)^Tx_i/{\hat
v}(\alpha \gamma(\hat \theta)^Tx_i)=0.
 \end{eqnarray}
Again we can also derive the similar results in Section~2.

In the above two
cases with given and unknown function $v(\cdot)$, the matrix $Q_1$, compared with its form in Theorem~\ref{theo1}, has a simpler structure because $v$ is a
function of $\beta_0$: $Q_1=:J(\theta_0)^T{\tilde
Q}J(\theta_0)=E\{g'(\beta_0^TX)^2/v(\beta_0^TX)J(\theta_0)^T
[X-E(X|\gamma_0^TX)][X-E(X|\gamma_0^TX)]^T J(\theta_0)\}.$

2. The limiting variances are related to re-parametrization. Thus, it is also  of interest to explore whether there is a re-parametrization so that the variances attain the minimum among all possible re-parametrization methods. It is noted that our results only need the orthogonality between $\gamma$ and its Jacobin matrix. Thus any re-parametrization with this property can be used in deriving similar results.

3. A more fundamental issue is about the asymptotic efficiency in the new framework we consider. For the GLM, we introduce a nuisance parameter and regard the link function as its true value. This is for the purpose of introducing automatically an extra ``estimation" of $g(\cdot)$ at the true value $\beta_0$. When this ``estimation" is non-positively correlated to  the one that can achieve the Fisher information, the resulting estimation can have chance to be more efficient.   This idea is general in principle. Hence, regardless of  re-parametrization, could we have an optimal selection of extra ``estimation" to achieve smallest possible variance? We guess that the variance we obtain by {\bf Procedure}~2 might be of some optimality property.  More generally, for an estimation of parameter of interest $\beta_0$, say, $T_n$, if we would be able to define an extra ``estimation"  $T_{n1}(\beta_0)$ at the true value of $\beta_0$ with non-positive correlation to $T_n-\beta_0$, the sum $(T_n-\beta_0)+T_{n1}(\beta_0)$ may have a smaller variance than that of $(T_n-\beta_0)$ when 2 time the covariance between them is, in absolute value, greater than the variance of $T_{n1}(\beta_0)$. As such, the key is how to find such an extra ``estimation". Our parameter space augmentation approach shows that it is possible at least for the GLM when we ``intentionally"  introduce a nuisance parameter into the model.

4. In this paper, we consider independent identically distributed cases. It may be readily extended to handle, say, cluster data and correlated data.

5. The advantage of the new method is its asymptotic efficiency.
Compared with  classical estimations, its computational cost should be more expensive. This is because it involves one-dimensional nonparametric
estimation and then we may face possible computational inefficiency and
instability with tuning parameter selection
when compared with the classical quasi-likelihood.
This deserves further studies
for finite sample implementation. 
\section{Appendix}\label{sec5}
\subsection{Conditions}
To obtain the asymptotic behavior of the estimators, we first give
the following conditions for technical purpose:
\newcounter{con}
\setcounter{con}{1}
\begin{list}{\bfseries \upshape C\arabic{con}.}
{\usecounter{con}} \item{(i)}\quad The distribution of $X$ has a
compact support set $A$.  \\(ii)\quad The density function
$f_\gamma(\cdot)$ of $\gamma^TX$  satisfies Lipschitz
condition of order $1$ for $\gamma$ in a neighborhood of $\gamma_0$.
Further, $\gamma_0^TX$ has a  bounded density function
$f_{\gamma_0}(\cdot)$ on its support $\mathcal{T}$.
\item{(i)}\quad
The function $g(\gamma^TX)=E(Y|\gamma^TX)$ has two bounded and
continuous derivatives.
\\(ii)\quad Let $l(z)=E(X/v(X)|\gamma_0^TX=z)$, and $l_s(\cdot)$ is the $s$-th component of $l(\cdot)$,  $1\leq s\leq p$. $l_s(\cdot)$ satisfies Lipschitz condition of order $1$. \\
(iii) \quad$Cov(X/v(X)|\gamma_0^TX=z)$ exists, $\max_{1\leq s\leq
p}E((X_s/v(X))^2|\gamma_0^TX=z)<c$.
\item{(i)}\quad The
kernel $K$ is a bounded, continuous and symmetric probability
density function, satisfying
\begin{eqnarray}
\int_{-\infty}^\infty u^2K(u)du\neq 0,\quad \int_{-\infty}^\infty
|u|^2K(u)du<\infty; \nonumber
\end{eqnarray}
\\(ii)\quad $K$ satisfies Lipschitz conditions on $\mathbf{R}^1$.
\item
$E(\varepsilon|X)=0,{\rm
var}(\varepsilon|X)=v(X)<\infty,E(\varepsilon^4|X)<\infty.$
\item $h\rightarrow0$, $nh^2/\log^2n\rightarrow\infty$,
$\limsup_{n\rightarrow\infty}nh^5\leq c<\infty$.
\item Both $V = E[g'^2(\beta_0^TX)XX^T/v(X)]$ and
$\tilde{V}=:  J(\theta_0)^TVJ(\theta_0) $ are positive definite
matrices, where $J(\theta_0)$ is the Jacobian matrix of
$\partial\gamma/\partial \theta$ evaluated at $\theta=\theta_0.$
\end{list}

\begin{rema} The Lipschitz condition and the derivatives in
C1 and C2 are standard smoothing conditions.  It is worth noticing that unlike other references, we do not bound
the density function of $X^T\beta_0$ from zero.  This is because when we use a truncation constant $c_0$, we can ensure the truncated density, and then
the denominators of $\hat{g}(z;\gamma_0,h)$ and
$\hat{g'}(z;\gamma_0,h)$ are bounded away from zero.
See
 Xia
and H\"ardle (2006) as reference.  C3 means that the kernel for
smoothing is of second order. C4 is a standard condition on the
error term. C6 ensures that the limiting variance for the estimator
$\gamma(\hat{\theta})$ exists.
\end{rema}

\subsection{Proofs}
\renewcommand{\theequation}{\mbox{A}.\arabic{equation}}
\setcounter{equation}{0}
As the following lemmas are very similar to those of Wang, Xue, Zhu
and Chong (2010), or Chang, Xue and Zhu (2010), we will only present
the main steps of proofs. A relevant reference is Zhu and Xue (2006). The main differences from existing ones in
these references are with different estimating equations and
estimations involved.

\begin{lemm}\label{lemma1} (Wang, Xue, Zhu and Chong 2010).Let
$\xi_1(x,\gamma),\ldots,\xi_n(x,\gamma)$ be a sequence of random
variables. Denote $f_{x,\gamma}(V_i)=\xi_i(x,\gamma)$ for
$i=1,\ldots,n$, where $V_1,\ldots,V_n$ be a sequence of random
variables, and $f_{x,\gamma}$ is a function on $\mathcal{A}_n$,
where $\mathcal{A}_n=\{(x,\gamma):(x,\gamma)\in A\times
R^p,\|\gamma-\gamma_0\|\leq cn^{-1/2}\}$ and $c>0$ is a constant.
Assume that $f_{x,\gamma}$ satisfies
\begin{eqnarray}\label{A1}
\frac{1}{n}\sum^n_{i=1}\left|
f_{x,\gamma}(V_i)-f_{x^\ast,\gamma^\ast}(V_i)\right|\leq
cn^\alpha\left(\parallel\gamma-\gamma^\ast\parallel+\parallel
x-x^\ast\parallel\right)
\end{eqnarray}
for some constants $x^\ast,\gamma^\ast$, $a>0$ and $c>0$. Let
$\varepsilon_n>0$ depend only on $n$. If
\begin{eqnarray}\label{A2}
P\left\{\left|\frac{1}{n}\sum^n_{i=1}\xi_i(x,\gamma)\right|>\frac{1}{2}\varepsilon_n\right\}\leq\frac{1}{2}
\end{eqnarray}
for $(x,\gamma)\in\mathcal{A}_n$, then we have
\begin{eqnarray}\label{A3}
&&P\left\{\sup_{(x,\gamma)\in\mathcal{A}_n}\left|\frac{1}{n}\sum^n_{i=1}\xi_i(x,\gamma)\right|>\frac{1}{2}\varepsilon_n\right\}\nonumber\\
&\leq&
c_1n^{2pa}\varepsilon^{-2p}_nE\left\{\sup_{(x,\gamma)\in\mathcal{A}_n}
2\exp(\frac{-n^2\varepsilon^{2}_n/128}{\sum^n_{i=1}\xi^2_i(x,\gamma)})\wedge
1\right\},
\end{eqnarray}
where $c_1>0$ is a constant.
\end{lemm}

\begin{lemm}\label{lemma2}
 Suppose that conditions {\rm C1,C2} and {\rm C3(i)} hold.
If $h=cn^{-a}$ for any $0<a<1/2$ and some constants $c>0$, then for
$i=1,\ldots,n$, we have
\begin{eqnarray}
&&E[g(\beta_0^Tx_j)-\sum_{i =1,i\neq j}^n
W_{ni}(\gamma_0^Tx_j;\gamma_0,h)g(\beta_0^Tx_i)]^2=O(h^4), \label{A4} \\
&&E[g(\alpha_0\gamma^Tx)-\sum_{i=1, i\neq j}^n
W_{ni}(\gamma^Tx;\gamma,h)g(\alpha_0\gamma^Tx_i)]^2=O(h^4),\\
&&E\left[\sum_{i=1}^n
W_{nj}(\gamma_0^Tx_i;\gamma_0)g'(\beta_0^Tx_i)l_s(\gamma_0^Tx_i)-
g'(\beta_0^Tx_j)l_s(\gamma_0^Tx_j)\right]^2=O(\sqrt{h})\nonumber,
\\
\end{eqnarray}
where $l_s$ is defined in Condition C2. 
\end{lemm}

{\it Proof}. The basic arguments of proof are from  Wang, Xue, Zhu
and Chong (2010), and it is very similar with Lemma 5.2 in Chang,
Xue and Zhu (2010). We then provide their main steps. We prove the
first conclusion and the others can be proven similarly. Denote
$z_i=\gamma_0^Tx_i$, $\tilde{z}_i=\beta_0^Tx_i$, replacing
$S_{n,r}(\gamma_0^Tx_i;\gamma_0,h)$ and
$U_{nj}(\gamma_0^Tx_i;\gamma_0,h)$ by $S_{n,r}(z_i)$ and
$U_{nj}(z_i)$ respectively for simplicity. We have
\begin{eqnarray}
g(\tilde{z}_j)-\sum_{i=1}^n
W_{ni}(z_j)g(\tilde{z}_i)=\frac{\sum_{i=1}^n
U_{ni}(z_j)[g(\tilde{z}_j)-g(\tilde{z}_i)]}{\sum_{i=1}^n
U_{ni}(z_j)}. \label{A8}
\end{eqnarray}

Let $E_{z_i}[\cdot]$ denote the conditional expectation given $z_i$.
Denote $T_n=O_r(a_n)$ if $E|T_n|^r=O(a^r_n)$. Using Cauchy-Schwarz
Inequality, we can easily obtain
\begin{eqnarray}
O_r(a_n)O_r(b_n)= O_{r/2}(a_nb_n),\label{A9}\\
T_n=E_{z_i}[T_n]+O_r((E|T_n-E_{z_i}[T_n]|^4)^{1/r}). \label{A10}
\end{eqnarray}

From condition C1(ii) we know
$M_{f_{\gamma_0}}=\sup_{t}f_{\gamma_0}(t)<\infty$, and there exists
a constant $L>0$ such that for any real $y$ and $t$,
$|f_{\gamma_0}(y)-f_{\gamma_0}(t)|<L|y-t|$. Using this facts. when
$n$ large enough, we can obtain that for $i=1,\ldots,n$,
\begin{eqnarray}\label{A11}
E_{z_i}[S_{n,r}(z_i)]=h^r\mu_rf_{\gamma_0}(z_i)(1+O(h)),\quad
r=0,1,2
\end{eqnarray}
where $\mu_r=\int^\infty_{-\infty}u^rK(u)du$, $r=0,1,2$, $\mu_0=1$
and $\mu_1=0$.

By the inequality of sum of independent random variables(see Petrov,
1995), we obtain that
\begin{eqnarray}
&&E\{|S_{n,r}(z_i)-E_{z_i}S_{n,r}(z_i)|^4\}=E\{E_{z_i}|S_{n,r}(z_i)-E_{z_i}S_{n,r}(z_i)|^4\}\nonumber\\
&\leq&
cn^{-4}(n-1)h^{4r-3}\int_{-\infty}^\infty|u^rK(u)|du({E}f_{\gamma_0}(z_i)+L{h})\nonumber\\
&&+cn^{-4}E[2(n-1)h^{2r-3}\int_{-\infty}^\infty u^{2r}K^2(u)du(f_{\gamma_0}(z_i)+L{h})]\nonumber\\
&=&O(n^{-3}h^{4r-3})+O(n^{-2}h^{2r-2}).
\end{eqnarray}
This together with (\ref{A10}) and (\ref{A11}) proves that
\begin{eqnarray}
S_{n,r}(z_i)&=&E_{z_i}[S_{n,r}(z_i)]+O_r(h^r(nh)^{-1/2})\nonumber\\
&=&\mu_rh^rf_{\gamma_0}(z_i)(1+O_4(h+(nh)^{-1/2})).\label{A13}
\end{eqnarray}
By (\ref{A9}) and (\ref{A13}), we have
\begin{eqnarray}
\frac{1}{n}\sum_{i=1}^n U_{ni}(z_j)&=&S_{n,0}(z_j)S_{n,2}(z_j)-S^2_{n,1}(z_j)\nonumber\\
&=&\mu_2h^2f^2_{\gamma_0}(z_j)(1+O_4(h+(nh)^{-1/2})).\label{A14}
\end{eqnarray}
Let $U_n(z)=(nh^2)^{-1}\sum_{i=1}^n U_{ni}(z)$ and
$U(z)=\mu_2f^2_{\gamma_0}(z)$. Using {Lemma \ref{lemma1}},
(\ref{A14}) and Borel-Cantelli's Lemma, we can prove[refer to Wang,
Xue, Zhu and Chong (2010)]
\begin{eqnarray}
\sup_{z\in\mathcal{T}}\left|U_n(z)-U(z)\right|\rightarrow 0,\quad
a.s.
\end{eqnarray}
From condition C1 we know $\inf_{z\in\mathcal{T}}f_{\gamma_0}(z)\geq
c>0$. thus, when $n$ is large enough,
\begin{eqnarray}
\inf_{z\in\mathcal{T}}|U_n(z)|\geq\inf_{z\in\mathcal{T}}|U(z)|-\sup_{z\in\mathcal{T}}|U_n(z)-U(z)|>\mu_2c^2/2>0,\quad
a.s. \label{A16}
\end{eqnarray}
Write
$H(\tilde{z}_i,\tilde{z}_j)=g(\tilde{z}_j)-g(\tilde{z}_i)+g'(\tilde{z}_j)(z_i-z_j)\alpha_0$.
Noting that $\sum_{i=1}^n U_{ni}(z_j)(z_i-z_j)=0$, we have
\begin{eqnarray}
&&\sum_{i=1}^n
U_{ni}(z_j)(g(\tilde{z}_j)-g(\tilde{z}_i))=\sum_{i=1}^n
H(\tilde{z}_i,\tilde{z}_j)U_{ni}(z_j)\nonumber\\
&=&\sum_{i=1}^n H(\tilde{z}_i,\tilde{z}_j)K_h(z_i-z_j)S_{n,2}(z_j)\nonumber\\
&&-\sum_{i=1}^n
H(\tilde{z}_i,\tilde{z}_j)(z_i-z_j)K_h(z_i-z_j)S_{n,1}(z_j).\label{A17}
\end{eqnarray}
Similar as the proof of (\ref{A13}) for $r=0,1$, we have
\begin{eqnarray}
&&\frac{1}{nh^{2+r}}\sum_{i=1}^n H(\tilde{z}_i,\tilde{z}_j)(z_i-z_j)^rK_h(z_i-z_j)\nonumber\\
&=&h^{-2-r}E_{z_j}[H(\tilde{z}_i,\tilde{z}_j)(z_1-z_j)^rK_h(z_1-z_j)]+o_4(1)\nonumber\\
&=:&d_{nr}+o_4(1),\quad 2\leq j\leq n,\label{A18}
\end{eqnarray}
and $|d_{nr}|\leq
c\int|u|^{2+r}K(u)du(f_{\gamma_0}(Z_j)+O(h))=O(1)$. This together
with (\ref{A13}), (\ref{A17}) and (\ref{A18}) proves
\begin{eqnarray}
\sum_{i=1}^n
U_{ni}(z_j)(g(\tilde{z}_j)-g(\tilde{z}_i))=nh^4\mu_2f_{\gamma_0}(Z_j)d_{n0}+o_2(nh^4).
\nonumber
\end{eqnarray}
Thus, combining again (\ref{A16}), we conclude
\begin{eqnarray}
&&E\left|\frac{\sum_{i=1}^n
U_{ni}(z_j)(g(\tilde{z}_j)-g(\tilde{z}_i))}{\sum_{i=1}^n
U_{ni}(z_j)}\right|^2\nonumber\\
&\leq&c(nh^2)^{-2}E\left|\sum_{i=1}^n
U_{ni}(z_j)(g(\tilde{z}_j)-g(\tilde{z}_i))\right|^2=O(h^4).
\end{eqnarray}
(\ref{A4}) now follows from (\ref{A8}) and (\ref{A18}). \hfill$\#$

\begin{lemm} \label{lemma3} Under the conditions of Lemma \ref{lemma2}, we have
\begin{eqnarray}\label{A20}
E\left\{\sum_{i=1}^n
W^2_{ni}(\gamma^Tx;\gamma,h)\right\}=O((nh)^{-1}),
\end{eqnarray}
\begin{eqnarray}
\left\{\begin{array}{l}
E\left\{W^2_{ni}(\gamma_0^Tx_i;\gamma_0,h)\right\}=O((nh)^{-2})\\
E\left\{\sum_{i=1,i\neq j}^n
W^2_{ni}(\gamma_0^Tx_j;\gamma_0,h)\right\}=O((nh)^{-1})\end{array}.\right.
\end{eqnarray}
\end{lemm}

{\it Proof}. The proof of Lemma \ref{lemma3} is similar as the proof
in Lemma \ref{lemma2}, hence, we omit it. \hfill$\#$

\begin{lemm}\label{lemma4}
Suppose that conditions {\rm C1-C4} and {\rm C5(i)} hold. We then
have
\begin{eqnarray}\label{A23}
\sup_{(x,\gamma)\in\mathcal{A}_n}|g(\beta_0^Tx)-\hat{g}(\gamma^Tx;\gamma,h)|=O_p((nh/\log
n)^{-1/2}).
\end{eqnarray}
\end{lemm}

{\it Proof}.  Write
$\tilde{g}(x_i,\varepsilon_i)=g(\beta_0^Tx)-g(\beta_0^Tx_i)-\varepsilon_i$,
$i=1,\ldots,n$. We have
\begin{eqnarray}
g(\beta_0^Tx)-\hat{g}(\gamma^Tx;\gamma,h)=\sum_{i=1}^n
W_{ni}(\gamma^Tx;\gamma,h)\tilde{g}(x_i,\varepsilon_i).
\end{eqnarray}
Let $\xi_i(x,\gamma)=n(nh/\log
n)^{1/2}W_{ni}(\gamma^Tx;\gamma,h)\tilde{g}(x_i,\varepsilon_i)$,
$f_{x,\gamma}(V_i)=\xi_i(x,\gamma)$, $V_i=(x_i,\varepsilon_i)$,
$i=1,\ldots,n$. Now we test and verify (\ref{A1}) and (\ref{A2}) in
Lemma \ref{lemma1}. A simple calculation yields (\ref{A1}). For
(\ref{A2}), by Lemma \ref{lemma2} and Lemma \ref{lemma3}(\ref{A20})
and noting that
$\sup_{(x,\gamma)\in\mathcal{A}_n}|g(\alpha_0\gamma^Tx)-g(x^T\beta_0)|=O(n^{-1/2})$,
$\gamma = \beta/\|\beta\|$, we have
\begin{eqnarray}
&&E[g(\beta_0^Tx)-\hat{g}(\gamma^Tx;\gamma,h)]^2\nonumber\\
& = &E[\sum_{i=1}^n
W_{ni}(\gamma^Tx;\gamma,h)\tilde{g}(x_i,\varepsilon_i)]^2\nonumber\\
&\leq& cE[g(\beta_0^Tx)-\sum_{i=1}^n
W_{ni}(\gamma^Tx;\gamma,h)g(\beta_0^Tx_i)]^2\nonumber\\
&& + {c}E[\sum_{i=1}^n
W^2_{ni}(\gamma^Tx;\gamma,h)] + O(n^{-1})\nonumber\\
&\leq&ch^4 + c(nh)^{-1}. \label{A26}
\end{eqnarray}
Given a $M>0$, by Chevbychev's inequality and (\ref{A26}), we have
\begin{eqnarray}
&&P\left\{\left|\frac{1}{n}\sum_{i=1}^n\xi_i(x,\gamma)\right|>\frac{1}{2}M\right\}\leq4M^{-2}E\left[\frac{1}{n}\sum_{i=1}^n
\xi_i(x,\gamma)\right]^2\nonumber\\
&\leq&4M^{-2}nh(\log n)^{-1}E\left[\sum_{i=1}^n W_{ni}(x^T\gamma;\gamma,h)\tilde{g}(X_i,\varepsilon_i)\right]^2\nonumber\\
&\leq&cM^{-2}(cnh^{5}+c(\log n)^{-1}).\label{A27}
\end{eqnarray}
Therefore, from condition C5, we can choose $M$ large enough so that the right
hand side of (\ref{A27}) is less than or equal to
$\displaystyle\frac{1}{2}$. Hence, (\ref{A2}) is satisfied. We now
can use (\ref{A3}) of Lemma \ref{lemma1} to get (\ref{A23}). By
(\ref{A20}), we obtain that
\begin{eqnarray}
n^{-2}\sum_{i=1}^n E\xi^2_i(x,\gamma)&= & nh(\log
n)^{-1}\sum_{i=1}^n
E[W_{ni}(\gamma^Tx;\gamma,h)\tilde{g}(x_i,\varepsilon_i)]^2\nonumber\\
& \leq & cnh(\log n)^{-1}\sum_{i=1}^n
EW^2_{ni}(\gamma^Tx;\gamma,h)\nonumber\\
& \leq & c(\log n)^{-1}. \nonumber
\end{eqnarray}
This implies that $n^{-2}\sum_{i=1}^n \xi^2_i(x,\gamma)=O_p((\log
n)^{-1})$. Hence, from Lemma \ref{lemma1} we have
\begin{eqnarray}
P\left\{\sup_{(x,\gamma)\in\mathcal{A}_n}\left|\frac{1}{n}\sum_{i=1}^n\xi_i(x,\gamma)\right|>\frac{1}{2}M\right\}\leq
cn^{2pa}M^{-2p}\exp(-cM^2\log n). \nonumber
\end{eqnarray}
The right hand side of the above formula to zero when $M$ is large
enough. There, (\ref{A23}) is shown. \hfill$\#$

\

\begin{lemm}\label{lemma5} Suppose that conditions {\rm C1-C6} are satisfied. Then we
have
\begin{eqnarray}
\sup_{\theta\in\mathcal{B}^\ast}\left\| R(\theta)-U(\theta_0)+
nV_1(\theta-\theta_0)\right\|=o_p(\sqrt{n}),
\end{eqnarray}
where $\mathcal{B}^\ast=\{\theta:\|\theta-\theta_0\|\leq
cn^{-1/2}\}$, for a constants $c>0$, $V_1$ is defined in condition
C6, and
\begin{eqnarray}
R(\theta)&=&\sum_{i=1}^n[y_i-\hat{g}(\gamma(\theta)^Tx_i;\theta,h)]g'(\hat{\alpha}\gamma(\theta)^Tx_i)J(\theta)^Tx_i/v(x_i)\\
U(\theta_0)&=&\sum_{i=1}^n\varepsilon_ig'(\beta_0^Tx_i)J(\theta_0)^T[x_i/v(x_i)-E(X/v(X)|\gamma_0^Tx_i)].
\end{eqnarray}
 as our denotation in section \ref{two-pahse}, $\beta_0 =
\alpha_0 \gamma_0$, $\gamma_0 = \gamma(\theta_0)$, and
$\hat{\alpha}$ is a consistent estimate of $\alpha_0$ in section
\ref{two-pahse}.
\end{lemm}

{\it Proof}. Separating $R(\theta)$, we have
\begin{eqnarray}\label{R}
R(\theta)&=&\sum_{i=1}^n\varepsilon_ig'(\beta_0^Tx_i)J^T(\theta)[x_i/v(x_i)-E(x_i/v(x_i)|\gamma_0^Tx_i)]\nonumber\\
&+&\sum_{i=1}^n\varepsilon_i[g'(\hat{\alpha}\gamma(\theta)^Tx_i)-g'(\alpha_0\gamma_0^Tx_i)]J^T(\theta)x_i/v(x_i)\nonumber\\
&-&\sum_{i=1}^n
g'(\beta_0^Tx_i)J^T(\theta)x_i/v(x_i)[\hat{g}(\gamma(\theta)^Tx_i;\theta,h)-\hat{g}(\gamma_0^Tx_i;\gamma_0,h)]\nonumber\\
&-&\sum_{i=1}^n
g'(\beta_0^Tx_i)J^T(\theta)\{x_i/v(x_i)[\hat{g}(\gamma_0^Tx_i;\gamma_0,h)-g(\beta_0^Tx_i)]-\varepsilon_il(\gamma_0^Tx_i)\}\nonumber\\
&-&\sum_{i=1}^n[\hat{g}(\gamma(\theta)^Tx_i;\theta,h)-g(\beta_0^Tx_i)][g'(\hat{\alpha}\gamma(\theta)^Tx_i)-
g'(\alpha_0\gamma_0^Tx_i)]J^T(\theta)x_i/v(x_i)\nonumber\\
&=:&R_1(\theta)+R_2(\theta)-R_3(\theta)-R_4(\theta)-R_5(\theta),
\end{eqnarray}
where $l(\gamma_0^TX)=E(X/V(X)|\gamma_0^TX)$ which is defined in C2.
$R_1(\theta)$ can be written as
\begin{eqnarray}
R_1(\theta) &=&
\sum_{i=1}^n\varepsilon_ig'(\alpha_0\gamma(\theta_0)^Tx_i)J(\theta_0)^T[x_i/v(x_i)-E(X/v(X)|\gamma_0^Tx_i)]\nonumber\\
&+&\sum_{i=1}^n\varepsilon_ig'(\alpha_0\gamma(\theta_0)^Tx_i)(J(\theta)-J(\theta_0))^T[x_i/v(x_i)-E(X/v(X)|\gamma_0^Tx_i)]\nonumber
\end{eqnarray}
Note that $\left\|J(\theta)-J(\theta_0)\right\|= O_p(n^{-1/2})$, for
all $\theta\in\mathcal{B}^\ast$, then by the law of large numbers,
we have
\begin{eqnarray}\label{R1}
\sup_{\theta\in\mathcal{B}^\ast}\left\|
R_1(\theta)-U(\theta_0)\right\|=o_p(\sqrt{n}).
\end{eqnarray}

For $R_2(\theta)$, as
$\left\|\gamma(\theta)-\gamma(\theta_0)\right\|= O_p(n^{-1/2})$,
$\left\|\hat{\alpha}\gamma(\theta) -
\alpha_0\gamma(\theta_0)\right\| = O_p(n^{-1/2})$, for all
$\theta\in\mathcal{B}^\ast$. By the Taylor expansion of
$g'(\hat{\alpha}\gamma(\theta)^Tx_i) -
g'(\alpha_0\gamma(\theta_0)^Tx_i)$, and the law of large numbers,
 we have
\begin{eqnarray}\label{R2}
\sup_{\theta\in\mathcal{B}^\ast}\left\|R_2(\theta)\right\| =
o_p(\sqrt{n})
\end{eqnarray}

By Taylor expansion of $\hat{g}$, we have for $R_3(\theta)$
\begin{eqnarray}
R_3(\theta)&=&\sum_{i=1}^n
g'(\alpha_0\gamma(\theta_0)^Tx_i)J(\theta)^Tx_i/v(x_i)[\hat{g}(\gamma(\theta)^Tx_i;\theta,h)
- \hat{g}(\gamma(\theta_0)^Tx_i;\theta_0,h)]\nonumber\\
&=& \sum_{i=1}^n
g'(\alpha_0\gamma(\theta_0)^Tx_i)J(\theta_0)^Tx_i/v(x_i)x_i^T[\hat{g}'(\gamma(\theta_0)^Tx_i;\theta_0,h)(\gamma(\theta)-\gamma(\theta_0))]
+ o_p(\sqrt{n})\nonumber
\end{eqnarray}
where $\hat{g}(\gamma(\theta_0)^Tx_i;\theta_0,h) =
\hat{E}(Y|\gamma(\theta_0)^Tx_i)$. Known that
$\hat{E}'(Y|\gamma(\theta_0)^Tx_i) \overset{p}{\rightarrow}
E'(Y|\gamma(\theta_0)^Tx_i)$, and actually
$E'(Y|\gamma(\theta_0)^Tx_i) = \frac{\partial
g(\beta_0^TX)}{\partial(\gamma(\theta_0)^TX)}|_{\gamma(\theta_0)^Tx_i}
= \alpha_0 g'(\alpha_0\gamma(\theta_0)^Tx_i)$. So, we obtain that
\begin{eqnarray}
R_3(\theta)&=& \alpha_0\sum_{i=1}^n
[g'(\beta_0^Tx_i)]^2J(\theta_0)^Tx_i/v(x_i)x_i^TJ(\theta_0)(\theta-\theta_0)+
o_p(\sqrt{n})\nonumber
\end{eqnarray}
by Weak Law of Large Numbers, we can derive
\begin{eqnarray}\label{R3}
\sup_{\theta\in\mathcal{B}^\ast}\left\|R_3(\theta)
-  n V_1(\theta -
\theta_0) \right\| = o_p(\sqrt{n})
\end{eqnarray}
where $V_1$ is defined in C6.

For $R_4(\theta)$, using the similar method as in Lemma A.5. of
Chang, Xue and Zhu(2010), using our Lemma \ref{lemma2} and Lemma
\ref{lemma3},  write $R_{4}(\theta)=J^T(\theta)R_{4}^\ast(\theta)$.
Let $R_{4s}^\ast$ denote the $s$th component of
$R_{4}^\ast(\theta)$, we can derive that
\begin{eqnarray}
\frac{1}{n}E\left({R_{4s}^\ast}^2\right)& =&
\frac{1}{n}E\left(\sum_{i=1}^n
g'(\beta_0^Tx_i)\{x_{is}/v(x_i)[\hat{g}(\gamma_0^Tx_i)-g(\beta_0^Tx_i)]-
\varepsilon_i l_s(\gamma_0^Tx_i)\} \right)^2\\\nonumber & \leq &
 cn^{-1}E\left(\sum_{i=1}^n g'(\beta_0^Tx_i)x_{is}/v(x_i)[\sum_{j=1}^n
 W_{nj}(\gamma_0^Tx_i)g(\beta_0^Tx_j)-g(\beta_0^Tx_i)]\right)^2\\\nonumber
 &+& cn^{-1}E\left( \sum_{i=1}^n g'(\beta_0^Tx_i)[x_{is}/v(x_i)\sum_{j=1}^n W_{nj}
 (\gamma_0^Tx_j)\varepsilon_j - \varepsilon_i l_s(\gamma_0^Tx_i)]
 \right)^2\\\nonumber
 &=:& cn^{-1}E (R_{4s1}^\ast)^2 + c n^{-1}E (R_{4s2}^\ast)^2
 \label{R42}
\end{eqnarray}
\begin{eqnarray}
E (R_{4s1}^\ast)^2 &\leq & c \sum_{i=1}^n\left( E
g'^2(\beta_0^Tx_i)[\sum_{j=1}^n W_{nj}(\gamma_0^Tx_i)g(\beta_0^Tx_j)
- g(\beta_0^Tx_i)]^2
E((x_{is}/v(x_i))^2|\gamma_0^Tx_i)\right)\nonumber
\end{eqnarray} by C(iii) and Lemma\ref{lemma2} , we have
\begin{eqnarray}\label{R4s1}
E (R_{4s1}^\ast)^2 &\leq & c n h^4
\end{eqnarray}
For $R_{4s2}^\ast$
\begin{eqnarray*}
E (R_{4s2}^\ast)^2 & = & \sum_{j=1}^n E\left( v(x_j)(\sum_{i=1}^n g'
(\beta_0^Tx_i)x_{is}/v(x_i)W_{nj}(\gamma_0^Tx_i) -
g'(\beta_0^Tx_j)l_s(\gamma_0^Tx_j))^2 \right)\\\nonumber &\leq& c
\sum_{i=1}^n E\left( v(x_i)(\sum_{j=1}^n
g'(\beta_0^Tx_j)W_{ni}(\gamma_0^Tx_j)l_s(\gamma_0^Tx_j) -
g'(\beta_0^Tx_i)l_s(\gamma_0^Tx_i))^2 \right)\\\nonumber &+& c
\sum_{i=1}^n E \left(v(x_i)(\sum_{j=1}^n
g'(\beta_0^Tx_j)W_{ni}(\gamma_0^Tx_j)(x_{js}/v(x_j)-l_s(\gamma^Tx_j)))^2
\right).
\end{eqnarray*}
BY  Lemma \ref{lemma2}, C(iii) and  Lemma \ref{lemma3}, we can
obtain
\begin{eqnarray}\label{R4s2}
E (R_{4s2}^\ast)^2 \leq c n \sqrt{h} + c n(nh)^{-1}
\end{eqnarray}
Together with (\ref{R42}), (\ref{R4s1}) and (\ref{R4s2}) we get
\begin{eqnarray*}
\frac{1}{n}E\left({R_{4s}^\ast}^2\right) \leq ch^4 + c \sqrt{h} +
c(nh)^{-1}\to 0.
\end{eqnarray*}

This implies
\begin{eqnarray}\label{R4}
\sup_{\theta\in\mathcal{B}^\ast}\left\|
R_{4}(\theta)\right\|=o_p(\sqrt{n}).
\end{eqnarray}

For $R_5(\theta)$, by Taylor expansion to $g'(\cdot)$, together with
$\left\|\hat{\alpha}\gamma(\theta) -
\alpha_0\gamma(\theta_0)\right\| = O_p(n^{-1/2})$, for all
$\theta\in\mathcal{B}^\ast$, we have
\begin{eqnarray}
R_5(\theta)&=& \sum_{i=1}^n
O_p(n^{-1/2})[\hat{g}(\gamma(\theta)^Tx_i;\theta,h)-g(\beta_0^Tx_i)]J^T(\theta)x_i/v(x_i).\nonumber
\end{eqnarray}
By Lemma \ref{lemma4}  and Condition C5, we obtain
\begin{eqnarray}\label{R5}
\sup_{\beta^{(r)}\in\mathcal{B}^\ast}\left\|
R_{5}(\beta^{(r)})\right\|=o_p(\sqrt{n}).
\end{eqnarray}
Together with (\ref{R1})- (\ref{R5}) we conclude the proof  of Lemma \ref{lemma5}.
\hfill$\#$

\begin{rema}
  When we consider (\ref{ngee1}) with $\tilde{I}_n(x_i)$, as
   $\|\tilde{I}_n(x_i)-I_n(x_i)\|\to 0$ can be proved, we can first
  substitute $\tilde{I}_n(x_i)$ with $I_n(x_i)$, and then use
  similar arguments as above, which will lead to  Remark
  \ref{rema1.1} right below Theorem \ref{theo1}.
\end{rema}

\ {\it Proof of Theorem~\ref{theo1}} The existence of the estimator
is easy to prove. As the estimating equation is  similar with the
one in Chang, Xue and Zhu(2010), and then the arguments are also very similar.
We omit the proof here.

\par Then we show the asymptotic normality.
\ $\gamma(\hat{\theta})$ is the solution of (\ref{ngee1}), that is,
$R(\hat{\theta})=0$. By Lemma \ref{lemma1}, we have
\begin{eqnarray}
R(\hat{\theta}) = U(\theta_0)- n V_1(\hat{\theta}-\theta_0) +
o_p(\sqrt{n}) = 0, \nonumber
\end{eqnarray}
and hence
\begin{eqnarray}
\sqrt{n}(\hat{\theta} - \theta_0) = V_1^{-1}n^{-1/2}U(\theta_0) +
o_p(1). \nonumber
\end{eqnarray}
It follows from (\ref{jacobin}) and (\ref{gamma1})  that
\begin{eqnarray}
\gamma(\hat{\theta})- \gamma(\theta_0) = J(\theta_0)(\hat{\theta} -
\theta_0) + o_p(n^{-1/2}). \nonumber
\end{eqnarray}
Thus, we have
\begin{eqnarray}
\sqrt{n}(\gamma(\hat{\theta}) - \gamma(\theta_0)) =
J(\theta_0)V_1^{-1}n^{-1/2}U(\theta_0) + o_p(1). \nonumber
\end{eqnarray}
Theorem \ref{theo1} now follows from Central Limit Theorems and
Slutsky's Theorem. \hfill$\Box$

{\it Proof of Proposition~\ref{prop}} By the definition of $\hat
\alpha_I$, it is easy to see that
\begin{eqnarray}
\sqrt n (\hat \alpha_I-\alpha_0)&=&\sqrt n ({\hat \beta}_{I}^T{\hat
\beta}_{I}-{ \beta}_{0}^T\beta_0)/(\|\hat \beta_{I}\|+\|
\beta_0\|)\nonumber\\
&=&\sqrt n(({\hat \beta}_{I}-{ \beta}_{0})^T\gamma_0+o(1)\nonumber
\end{eqnarray}
Then it is clear that the limiting variance is
$\gamma_0^TV^{-1}\gamma_0$. For $\hat \gamma={\hat \beta}_I/\|{\hat
\beta}_I\|$, some elementary calculation yields that
\begin{eqnarray}
\hat \gamma-\gamma_0&=&({\hat \beta}_{I}/\|\hat \beta_{I}\|-{
\beta}_{0}/\| \beta_0\|)=({\hat \beta}_{I}\|\beta_{0}\|-{
\beta}_{0}\|\hat \beta_I\|)/(\|\hat \beta_{I}\|\|{
\beta}_{0} \|) \nonumber\\
&=&((I_p-\gamma_0\gamma_0^T)({\hat \beta}_{I}-{ \beta}_{0}))/\|{
\beta}_{0} \|+o(1/\sqrt n).\nonumber
\end{eqnarray}
From this presentation, we have the limiting variance:
$$
(I_p-\gamma_0\gamma_0^T)V^{-1}(I_p-\gamma_0\gamma_0^T)/\alpha_0^2.
$$
Note that
$V^{-1}=A(\theta_0)(A(\theta_0)^TVA(\theta_0))^{-1}A(\theta_0)^T$,
and $(I_p-\gamma_0\gamma_0^T)A(\theta_0)=({\bf 0}_p, J(\theta_0))$.
The result follows. \hfill$\Box$
\

{\it Proof of Theorem~\ref{theo2}} As both $\hat \alpha_I$ and
$\gamma(\hat \theta)$ are respectively convergent to $\alpha_0$ and
$\gamma_0$, it is easy to see that $\hat \alpha_I \gamma(\hat
\theta)-\beta_0=(\hat
\alpha_I-\alpha_0)\gamma_0+\alpha_0(\gamma(\hat
\theta)-\gamma_0)+o_p(1/\sqrt n)$. From the standard proof for the
asymptotic representation of $\hat \beta_{I}$ by the quasi
likelihood, we have the following result. Noting that $\hat \alpha
=\|\hat \beta_{I}\|$ and $\alpha_0 =\| \beta_0\|$, and together with
the proof of Proposition~\ref{prop}, the estimator by the
quasi-likelihood satisfies
\begin{eqnarray}
\sqrt n (\|\hat \beta_{I}\|-\| \beta_0\|)&=&\gamma_0^T\sqrt n(\hat \beta_{I}-\beta_0)+o_p(1)\nonumber\\
&=&\gamma_0^TV^{-1}\frac1{\sqrt n}\sum_{i=1}^n\varepsilon_ig'(\beta_0^Tx_i)x_i/v(x_i)+o_p(1).
\end{eqnarray}
Then the application of Central Limit Theorems yields that it converges in distribution to a normal distribution with mean zero and variance $\gamma_0^TV^{-1}\gamma_0$.
Further, from  the proof of Theorem~\ref{theo1}, we have
\begin{eqnarray}\label{A38}
&&\sqrt n (\gamma(\hat \theta)-\gamma_0)\nonumber\\
&=&J(\theta_0)V_1^{-1}\frac1{\sqrt
n}\sum_{i=1}^n\varepsilon_ig'(\beta_0^Tx_i)
J(\theta_0)^T[x_i/v(x_i)-E(X/v(X)|\gamma_0^Tx_i)]\nonumber\\
&&+o_p(1).
\end{eqnarray}
Note that $A(\theta_0)^T\gamma_0=(1, {\bf 0}_{p-1}^T)^T$ and
$A(\theta_0)^TJ(\theta_0)=({\bf 0}_{p-1},
J(\theta_0)^TJ(\theta_0))^T.$ Thus, our estimator $\hat{\beta} =
\hat{\alpha} \gamma(\hat{\theta})$ has
\begin{eqnarray}
&&\sqrt n  (\hat \beta-\beta_0)=\sqrt n (A(\theta_0)^T)^{-1}A(\theta_0)^T (\hat \beta-\beta_0)\nonumber\\
&=&(A(\theta_0)^T)^{-1}\frac1{\sqrt
n}\sum_{i=1}^n\varepsilon_ig'(\beta_0^Tx_i)\Big ((1,
{\bf 0}_{p-1}^T)^T\gamma_0^TV^{-1}x_i/v(x_i)+\nonumber\\
&& \qquad \qquad \qquad ({\bf 0}_{p-1}, J(\theta_0)^TJ(\theta_0))^T
\tilde{V}^{-1} J(\theta_0)^T[x_i/v(x_i)-E(X/v(X)|
\gamma_0^Tx_i)]\Big )\nonumber\\
&&+o_p(1).
\end{eqnarray}
Central Limit Theorems imply the asymptotic normality  with mean
zero and variance-covariance matrix
\begin{eqnarray*}
&&(A(\theta_0)^T)^{-1}\left(
        \begin{array}{cc}
          \gamma_0^TV^{-1}\gamma_0 & {\bf 0}_{p-1}^T \\
          {\bf 0}_{p-1} & J(\theta_0)^TJ(\theta_0)\tilde{V}^{-1}Q_1\tilde{V}^{-1}J(\theta_0)^TJ(\theta_0)\\
        \end{array}
      \right) A(\theta)^{-1}\\
&=&(A(\theta_0)^T)^{-1}\left (
\begin{array}{cc}
1 & {\bf 0}_{p-1}^T \\
{\bf 0}_{p-1} & J(\theta_0)^TJ(\theta_0)
\end{array}
      \right)
\left(
        \begin{array}{cc}
          \gamma_0^TV^{-1}\gamma_0 & {\bf 0}_{p-1}^T \\
          {\bf 0}_{p-1} & \tilde{V}^{-1}Q_1\tilde{V}^{-1}\\
        \end{array}
      \right) \left (
\begin{array}{cc}
1 & {\bf 0}_{p-1}^T \\
{\bf 0}_{p-1} & J(\theta_0)^TJ(\theta_0)
\end{array}
      \right)A(\theta_0)^{-1}\\
      &&=A(\theta_0)W A(\theta_0)^{T}.
\end{eqnarray*}
The last equation holds because $(A(\theta_0)^T)^{-1}\left (
\begin{array}{cc}
1 & {\bf 0}_{p-1}^T \\
{\bf 0}_{p-1} & J(\theta_0))^TJ(\theta_0))
\end{array}
      \right)=(A(\theta_0)^T)^{-1}A(\theta_0)^TA(\theta_0).$
The conclusion is proved. \hfill $\fbox{}$

{\it Proof of Theorem~\ref{theo3}}. From the proof of Theorem~\ref{theo2}, we only need to prove that $(\gamma_0^TV\gamma_0)^{-1}+
          \frac {\gamma_0^TVJ(\theta_0)\tilde{V}^{-1}Q_1\tilde{V}^{-1}J(\theta_0)^TV\gamma_0(\theta_0)}
          {(\gamma_0^TV\gamma_0)^2} $ is  the limiting variance of the estimator $\hat \alpha_F$.
\ From (\ref{gee2}), we can derive that
\begin{eqnarray}
&&\frac1{\sqrt n}\sum_{i=1}^n\varepsilon_i g'(  \alpha_0 \gamma_0^Tx_i) \gamma_0^Tx_i/v(x_i)\nonumber\\
&=&\gamma_0^TV\sqrt n(\hat \alpha_F \gamma(\hat \theta)-\beta_0)+o_p(1)\nonumber\\
&=&\gamma_0^TV\gamma_0\sqrt n(\hat \alpha_F-\alpha_0)
+\alpha_0\gamma_0^TV\sqrt n(\gamma(\hat
\theta)-\gamma_0)+o_p(1)\nonumber
\end{eqnarray}
Together with (\ref{A38}), we have
\begin{eqnarray}
&& \gamma_0^TV\gamma_0\sqrt n(\hat \alpha_F-\alpha_0)\nonumber\\
&=&\frac1{\sqrt n}\sum_{i=1}^n\varepsilon_i g'(  \beta_0^Tx_i) \gamma( \theta_0)^Tx_i/v(x_i)\nonumber\\
&-&\gamma_0^TVJ(\theta_0) \tilde{V}^{-1}\frac1{\sqrt
n}\sum_{i=1}^n\varepsilon_ig'(\beta_0^Tx_i)
J(\theta_0)^T[x_i/v(x_i)-E(X/v(X)|\gamma_0^Tx_i]\nonumber\\
\nonumber
\end{eqnarray}
The right hand side is of a limiting variance as
$$
\gamma_0^TV\gamma_0+\gamma_0^TVJ(\theta_0)
\tilde{V}^{-1}Q_1\tilde{V}^{-1}J(\theta_0)^TV\gamma_0-2\gamma_0^TVJ(\theta_0)
\tilde{V}^{-1}J(\theta_0)^TQ_2\gamma_0,
$$
 is the limiting variance-covariance matrix of $\gamma_0^TV\gamma_0\sqrt n(\hat
 \alpha_F-\alpha_0)$, where $Q_2 = E\left(g'^2(\beta_0^TX)(X/v(X)- E(X/v(X)|\gamma_0^TX))X^T
 \right)$. It is easy to prove that \newline$C  \gamma_0^T VJ(\theta_0)\tilde{V}J(\theta_0)^TQ_2\gamma_0 =
 0$. The proof is completed.
\hfill $\fbox{}$

{\it Proof of Theorem~\ref{theo4}}. 
From the structure of $V$ it is easy to see that $V_2$ can be
written as
$$V_2=A(\theta_0)^TVA(\theta_0)=\left(
        \begin{array}{cc}
          \gamma_0^TV\gamma_0&
          { \gamma_0^TVJ(\theta_0)} \\
          {J(\theta_0)^TV\gamma_0} & J(\theta_0)^TVJ(\theta_0) \\
        \end{array}
      \right)=:
\left(
        \begin{array}{cc}
          V_{11}&
          V_{12} \\
         V_{21} & \tilde{V} \\
        \end{array}
      \right).
      $$
      Further, we have
$$V_2^{-1}=\left(
        \begin{array}{cc}
          V^{-1}_{11}+V_{11}^{-2}V_{12}V_{22\cdot 1}^{-1}V_{21}&
          -V_{11}^{-1}V_{12}V_{22\cdot 1}^{-1} \\
         -V_{22\cdot 1}^{-1} V_{21}V_{11}^{-1} & V_{22\cdot 1}^{-1} \\
        \end{array}
      \right)
      $$
where $V_{22\cdot 1}=\tilde{V}-V_{21}V_{11}^{-1}V_{12}$. Note that
$V_{21}V_{11}^{-1}V_{12}$ is a non-negative semidefinite and
$\tilde{V}$ is positive definite. Thus, $V_{22\cdot 1}\le \tilde{V}$
and $V_{22\cdot 1}^{-1}\ge \tilde{V}^{-1}\ge
\tilde{V}^{-1}Q_1\tilde{V}^{-1}$ as $Q_1\le \tilde{V}$. Note that
$V^{-1}=A(\theta_0)V_2^{-1}A(\theta_0)^T$, and then $\gamma_0^T
V^{-1}\gamma_0=(V_2^{-1})_{11}=V^{-1}_{11}+V_{11}^{-2}V_{12}V_{22\cdot
1}^{-1}V_{21}$. These result in that $trace(V_2^{-1})\ge trace( W
)$. We then derive that, noting that
$A(\theta_0)^TA(\theta_0)=\left(
        \begin{array}{cc}
          1&
          {\bf 0}_{p-1}^T \\
         {\bf 0}_{p-1} & J(\theta_0)^TJ(\theta_0) \\
        \end{array}
      \right)$ is a diagonal matrix,
\begin{eqnarray*}trace(V^{-1})&=&
trace(A(\theta_0)V_2^{-1}A(\theta_0)^T)\\
&=&trace(V_2^{-1}A(\theta_0)^TA(\theta_0))\\
&=&V^{-1}_{11}+V_{11}^{-2}V_{12}V_{22\cdot
1}^{-1}V_{21}+trace(V_{22\cdot 1}^{-1}J(\theta_0)^TJ(\theta_0))\\
&\ge &
\gamma_0^TV^{-1}\gamma_0+trace(\tilde{V}^{-1}Q_1\tilde{V}^{-1}J(\theta_0)^TJ(\theta_0))\\
&=&trace(A(\theta_0)WA(\theta_0)^T).
 \end{eqnarray*}

 For $W_1$ we have
\begin{eqnarray*}
&&(\gamma_0^TV\gamma_0)^{-1}+
          \frac {\gamma_0^TVJ(\theta_0)\tilde{V}^{-1}Q_1\tilde{V}^{-1}J(\theta_0)^TV\gamma_0}
          {(\gamma_0^TV\gamma_0)^2}\\
&=&V_{11}^{-1}+V_{12}\tilde{V}^{-1}Q_1\tilde{V}^{-1}V_{21}V_{11}^{-2}\le
V_{11}^{-1}+V_{12}\tilde{V}^{-1}V_{21}V_{11}^{-2}\le
V_{11}^{-1}+V_{12}V_{22\cdot 1}^{-1}V_{21}V_{11}^{-2}.
\end{eqnarray*}
These result in that $trace(V_2^{-1})\ge trace( W_1 )$ and then
$$
trace(V^{-1})\ge trace(A(\theta_0)W_1A(\theta_0)^T).
$$
The proof is complete. \hfill $\fbox{}$

{}

\end{document}